\documentclass[a4paper,11pt]{article}
\usepackage{amsmath,amsthm,amssymb,amsfonts,fullpage,graphicx,pstricks}

\newtheorem{thm}{Theorem}[section]

\newtheorem{cor}[thm]{Corollary}

\newtheorem{lem}[thm]{Lemma}

\newtheorem{rem}[thm]{Remark}

\numberwithin{equation}{section}

\begin{document}

\title{Biased random walk\\on critical Galton-Watson trees\\conditioned to survive}
\author{D. A. Croydon\footnote{Dept of Statistics, University of Warwick, Coventry, CV4 7AL, United Kingdom; d.a.croydon@warwick.ac.uk.},
~A. Fribergh\footnote{CIMS, 251 Mercer Street, New York University, New York 10012-1185, U.S.A.; fribergh@cims.nyu.edu.} ~and ~T. Kumagai\footnote{RIMS, Kyoto University, Kyoto 606-8502, Japan; kumagai@kurims.kyoto-u.ac.jp.}}
\maketitle

\begin{abstract} We consider the biased random walk on a critical Galton-Watson tree conditioned to survive, and confirm that this model with trapping belongs to the same universality class as certain one-dimensional trapping models with slowly-varying tails. Indeed, in each of these two settings, we establish closely-related functional limit theorems involving an extremal process and also demonstrate extremal aging occurs.
\end{abstract}

\section{Introduction}

Biased random walks in inhomogeneous environments are a natural setting to witness trapping phenomena. In the case of supercritical Galton-Watson trees with leaves (see \cite{AFGH}, \cite{BH}, \cite{LPP}) or the supercritical percolation cluster on $\mathbb{Z}^d$ (see \cite{FH}), for example, it has been observed that dead-ends found in the environment can, for suitably strong biases, create a sub-ballistic regime that is characteristic of trapping. More specifically, for both of these models, the distribution of the time spent in individual traps has polynomial tail decay, and this places them in the same universality class as the one-dimensional heavy-tailed trapping models considered in \cite{Zindy}. Indeed, although the case of a deterministically biased random walk on a Galton-Watson tree with leaves is slightly complicated by a certain lattice effect, which means it can not be rescaled properly \cite{AFGH}, in the case of randomly biased random walks on such structures, it was shown in \cite{BH} that precisely the same limiting behaviour as the one-dimensional models of \cite{ESZ2} and \cite{Zindy} occurs. Moreover, there is evidence presented in \cite{FH} that suggests the biased random walk on a supercritical percolation cluster also has the same limiting behaviour. The universality class that connects these models was previously investigated in \cite{BC}, \cite{bacone} and \cite{bactwo}, and is characterised by limiting stable subordinators and aging properties.

The aim of this paper is to investigate biased random walks on critical structures. To this end, we choose to study the biased random walk on a  critical Galton-Watson tree conditioned to survive. With the underlying environment having radically different properties from its supercritical counterpart, we would expect different limiting behaviour, with more extreme trapping phenomena, to arise. It is further natural to believe that some of the properties of the biased random walk on the incipient infinite cluster for critical percolation on $\mathbb{Z}^d$, at least in high dimensions, would be similar to the ones proved in our context, as is observed to be the case for the unbiased random walk (compare, for instance, the results of \cite{BarKum} and \cite{KozNac}). Nevertheless, our current understanding of the geometry of this object is not sufficient to extend our results easily, and so we do not pursue this inquiry here. In particular, we anticipate that, as indicated by physicists in \cite{BD}, for percolation close to criticality there is likely to be an additional trapping mechanism that occurs due to spatial considerations, which means that, even without taking the effect of dead-ends into account, it is more likely for the biased random walk to be found in certain regions of individual paths than others (see \cite{BRWRRW} for a preliminary study in this direction).

Our main model -- the biased random walk on critical Galton-Watson trees conditioned to survive -- is presented in the next section, along with a summary of the results we are able to prove for it. This is followed in Section \ref{onedsec} with an introduction to a one-dimensional trapping model in which the trapping time distributions have slowly-varying tails. This latter model, which is of interest in its own right, is of particular relevance for us, as it allows us to comprehensively characterise the universality class into which the Galton-Watson trees we consider fall. Furthermore, the arguments we apply for the one-dimensional model provide a useful template for the more complicated tree framework.

\subsection{Biased random walk on critical Galton-Watson trees}
\label{treeint}

Before presenting the Galton-Watson tree framework, we recall some classical results for sums of random variables whose distribution has a slowly-varying tail. Let $(X_i)_{i=1}^\infty$ be independent random variables, with distributional tail $\bar{F}(u)=1-F(u)=\mathbf{P}(X_i>u)$ satisfying: $\bar{F}(0)=1$, $\bar{F}(u)>0$ for all $u>0$,
\begin{equation}\label{svartail}
\lim_{u\rightarrow\infty}\frac{\bar{F}(uv)}{\bar{F}(u)}=1,
\end{equation}
for any $v>0$, and $\bar{F}(u)\rightarrow 0$ as $u\rightarrow\infty$. A typical example is when the distribution in question decays logarithmically slowly, such as
\begin{equation}\label{taileg}
\bar{F}(u)\sim \frac{1}{(\ln u)^{\gamma}},
\end{equation}
for some $\gamma>0$, where throughout the article $f\sim g$ will mean $f(x)/g(x)\rightarrow 1$ as $x\rightarrow\infty$. A first scaling result for sums of the form $\sum_{i=1}^nX_i$ was obtained in \cite{Darling}, and this was subsequently extended by \cite{Kasa} to a functional result. In particular, in \cite{Kasa} it was established that if $L(x):=1/\bar{F}(x)$, then
\begin{equation}\label{mlim}
\left(\frac{1}{n}L\left(\sum_{i=1}^{nt}X_i\right)\right)_{t\geq 0}\rightarrow \left(m(t)\right)_{t\geq 0}
\end{equation}
in distribution with respect to the Skorohod $J_1$ topology (as an aid to the reader, we provide in the appendix a definition of the Skorohod $J_1$ and
$M_1$ topologies, the latter of which is applied in several subsequent results), where  $m=(m(t))_{t\geq 0}$ is an extremal process. To define $m$ more precisely, suppose that $(\xi(t))_{t\ge 0}$ is the symmetric Cauchy process, i.e., the L\'evy process with L\'evy measure given by $\mu((x,\infty))=x^{-1/2}$ for $x>0$, and then set
\[m(t)=\max_{0<s\le t}\Delta\xi(s),\]
where $\Delta\xi(s)=\xi(s)-\xi(s^-)$. (Observe that $(m(t))_{t\ge 0}$ is thus the maximum process of the Poisson point process with intensity measure $x^{-2}dxdt$.) We will prove that, in addition to appearing in the limit at (\ref{mlim}), this extremal process arises in the scaling limits of a biased random walk on a critical Galton-Watson tree and, as is described in the next section, a one-dimensional directed trap model whose holding times have a slowly-varying mean.

We continue by introducing some relevant branching process and random walk notation, following the presentation of \cite{CroyKum}. Let $Z$ be a critical ($\mathbf{E}Z=1$) offspring distribution in the domain of attraction of a stable law with index $\alpha\in (1,2]$, by which we mean that there exists a sequence $a_n\uparrow \infty$ such that
\begin{equation}\label{stabledomain}
\frac{Z[n]-n}{a_n}\buildrel{d}\over{\rightarrow} X,
\end{equation}
where $Z[n]$ is the sum of $n$ i.i.d. copies of $Z$ and $\mathbf{E}(e^{-\lambda X})=e^{-\lambda^\alpha}$ for $\lambda\geq 0$. Note that, by results of \cite[Chapters XIII and XVII]{Feller2}, this is equivalent to the probability generating function of $Z$ satisfying
\begin{equation}\label{genfun}
f(s):=\mathbf{E}(s^Z)=\sum_{k=0}^\infty p_ks^k=s+(1-s)^{\alpha}L(1-s),\hspace{20pt}\forall s\in (0,1),
\end{equation}
where $L(x)$ is slowly varying as $x\rightarrow 0^+$, and the non-triviality condition $\mathbf{P}(Z=1)\neq 1$ holding. We point out that the condition $\mathbf{E}(Z^2)<\infty$ is sufficient for the previous statements to hold with $\alpha=2$.

Denote by $(Z_n)_{n\geq 0}$ the corresponding Galton-Watson process, started from $Z_0=1$. It has been established in \cite[Lemma 2]{Slack} that if $q_n:=\mathbf{P}(Z_n>0)$, then
\begin{equation}\label{probdecay}
q_n^{\alpha-1}L(q_n)\sim \frac{1}{(\alpha-1)n},
\end{equation}
as $n\rightarrow\infty$, where $L$ is the function appearing in (\ref{genfun}). It is also well known that the branching process $(Z_n)_{n\geq 0}$ can be obtained as the generation size process of a Galton-Watson tree, $\mathcal{T}$ say, with offspring distribution $Z$. In particular, to construct the random rooted graph tree $\mathcal{T}$, start with a single ancestor (or root), and then suppose that individuals in a given generation have offspring independently of the past and each other according to the distribution of $Z$, see \cite[Section 3]{rrt} for details. The vertex set of $\mathcal{T}$ is the entire collection of individuals, edges are the parent-offspring bonds, and $Z_n$ is the number of individuals in the $n$th generation of $\mathcal{T}$. From (\ref{probdecay}), it is clear that $\mathcal{T}$ will be a finite graph $\mathbf{P}$-a.s. However, in \cite{Kesten}, Kesten showed that it is possible to make sense of conditioning $\mathcal{T}$ to survive or `grow to infinity'. More specifically, there exists a unique (in law) random infinite rooted locally-finite graph tree $\mathcal{T}^*$ that satisfies, for any $n\in\mathbb{Z}_+$,
\[\mathbf{E}\left(\phi(\mathcal{T}^*|_n)\right)=\lim_{m\rightarrow\infty}\mathbf{E}\left(\phi(\mathcal{T}|_n)|Z_{m+n}>0\right),\]
where $\phi$ is a bounded function on finite rooted graph trees of $n$ generations, and $\mathcal{T}|_n$, $\mathcal{T}^*|_n$ are the first $n$ generations of $\mathcal{T}$, $\mathcal{T}^*$ respectively. We will write $d_{\mathcal{T}^*}$ to represent the shortest path graph distance on $\mathcal{T}^*$.

Given a particular realisation of $\mathcal{T}^*$, we will denote by $X=((X_n)_{n\geq 0}, P_x^{\mathcal{T}^*},x\in \mathcal{T}^*)$ the discrete-time biased random walk on $\mathcal{T}^*$, and define this as follows. First, fix a bias parameter $\beta>1$, and assign to each edge connecting a vertex $x$ in generation $k$ to a vertex $y$ in generation $k+1$ a conductance $c(x,y):=\beta^k=:c(y,x)$. The transition probabilities of $X$ are then determined by
\[P^{\mathcal{T}^*}(x,y):=\frac{c(x,y)}{\sum_{y'\sim x}c(x,y')},\hspace{20pt}\forall x\sim y,\]
where the notation $x\sim y$  means that $x$ and $y$ are connected by an edge in $\mathcal{T}^*$. Thus, when at a vertex $x$ that is not equal to the root of $\mathcal{T}^*$, the probability of jumping to a neighbouring vertex further away from the root than $x$ is  $\beta$ times more likely than jumping towards the root. Using the usual terminology for random walks in random environments, we will say that $P_x^{\mathcal{T}^*}$ is the quenched law of the biased random walk on $\mathcal{T}^*$ started from $x$. Moreover, we introduce the annealed law for the process started from $\rho$, the root of the tree $\mathcal{T}^*$, by setting
\begin{equation}\label{annealedlaw}
\mathbb{P}_\rho(\cdot):=\int P_\rho^{\mathcal{T}^*}(\cdot) {\rm d}\mathbf{P}.
\end{equation}
It will be this law under which we investigate the rate at which the process $X$, which we call the biased random walk on a critical Galton-Watson tree conditioned to survive, escapes from the root.

The main result we prove for the process $X$ concerns the time it takes to progress along the backbone. To be more specific, as is described in more detail in Section \ref{strucsec}, $\mathbf{P}$-a.s. the tree $\mathcal{T}^*$ admits a unique backbone, that is, a semi-infinite path starting from the root, $\{\rho=\rho_0,\rho_1,\rho_2,\dots\}$ say. We define $(\Delta_n)_{n\geq 0}$ by setting
 \begin{equation}\label{xhit}
\Delta_n:=\inf\left\{m\geq 0:\:X_m=\rho_n\right\}
\end{equation}
to be the first time the process $X$ reaches level $n$ along this path. For this process, we are able to prove the following functional limit theorem.

{\thm \label{dngrowth} Let $\alpha\in(1,2]$. As $n\rightarrow\infty$, the laws of the processes
\[\left(\frac{(\alpha-1)\ln_+ \Delta_{nt}}{n\ln \beta}\right)_{t\geq 0}\]
under $\mathbb{P}_\rho$ converge weakly with respect to the Skorohod $J_1$ topology on $D([0,\infty),\mathbb{R})$ to the law of $(m(t))_{t\geq0}$.}
\bigskip

It is interesting to observe that this result is extremely explicit compared to its supercritical counterparts. Indeed, notwithstanding the fact the lattice-effect that was the source of somewhat complicated behaviour in \cite{AFGH} does not occur in the critical setting, the above scaling limit  clearly describes the $\beta$-dependence of the relevant slowdown effect. Note that, unlike in the supercritical case where there is a ballistic phase, this slowdown effect occurs for any non-trivial bias parameter, i.e. for any $\beta>1$. Furthermore, we remark that the dependence on $\alpha$ is natural: as $\alpha$ decreases and the leaves get thicker (in the sense that tree's Hausdorff dimension of $\alpha/(\alpha-1)$ increases, see \cite{LegallDuquesne, hmg}), the biased random walk moves more slowly away from its start point.

As suggested by comparing Theorem \ref{dngrowth} with (\ref{mlim}), the critical Galton-Watson tree case is closely linked with a sum of independent and identically-distributed random variables where $\bar{F}$ is  asymptotically equivalent to $\ln \beta/(\alpha-1)\ln x$. Although the logarithmic rate of decay is relatively easy to guess, finding the correct constant is slightly subtle, particularly for $\alpha\neq 2$. This is because, unlike in the supercritical case and the critical case with $\alpha=2$, when $\alpha\neq 2$ it can happen that there are multiple deep traps emanating from a single backbone vertex. As a result, we have to take special care which of these have actually been visited when determining the time spent there, meaning that the random variable which actually has the $\ln \beta/(\alpha-1)\ln x$ tail behaviour is not environment measurable (see Lemma \ref{ttail}). To highlight the importance of this consideration, which is also relevant albeit in a simpler way for $\alpha=2$, in Theorem \ref{edngrowth} we show that the constant that appears differs by a factor $\alpha$ when $\Delta_n$ is replaced by its quenched mean $E^{\mathcal{T}^*}_\rho\Delta_n$.

Theorem \ref{dngrowth} readily implies the following corollary for the projection, $(\pi(X_m))_{m\geq 0}$, of the process $(X_m)_{m\geq 0}$ onto the backbone (roughly, $\pi(X_m)$ is the vertex on the backbone from which the trap $X_m$ is located in emanates, see Section \ref{rwsec} for a precise definition). To state this, we define the right-continuous inverse $(m^{-1}(t))_{t\geq0}$ of $(m(t))_{t\geq 0}$ by setting
\begin{equation}\label{minv}
m^{-1}(t):=\inf\left\{s\geq 0:\:m(s)>t\right\}.
\end{equation}

{\cor \label{treecor} Let $\alpha\in(1,2]$. As $n\rightarrow\infty$, the laws of the processes
\[\left(\frac{ d_{\mathcal{T}^*}\left(\rho,\pi(X_{e^{nt}})\right)\ln \beta}{(\alpha-1) n}\right)_{t\geq 0}\]
under $\mathbb{P}_\rho$ converge weakly with respect to the Skorohod $M_1$ topology on $D([0,\infty),\mathbb{R})$ to the law of $(m^{-1}(t))_{t\geq0}$.}

{\rem Since the height of the leaves in which the random walk can be found at time $e^{n}$ (see the localisation result of Lemma \ref{localtree}) will typically be of order $n$, some further argument will be necessary to deduce a limit result for the graph distance $d_{\mathcal{T}^*}(\rho,X_n)$ itself.}
\bigskip

Another characteristic property that we are able to show is that the random walk also exhibits extremal aging.

{\thm \label{extagingtree} Let $\alpha\in(1,2]$. For any $0<a<b$, we have
\[\lim_{n\rightarrow\infty}\mathbb{P}_\rho\left(\pi(X_{e^{an}})=\pi(X_{e^{bn}})\right)=\frac{a}{b}.\]}

Although regular aging has previously been observed for random walks in random environments in the sub-ballistic regime on $\mathbb{Z}$ (see \cite{esz}), as far as we know, this is the first example of a random walk in random environment where extremal aging has been proved. As already hinted at, this kind of behaviour, as well as that demonstrated in Theorem \ref{dngrowth} and Corollary \ref{treecor}, places the biased random walk on a critical Galton-Watson tree conditioned to survive in a different universality class to that of the supercritical structures discussed previously. In the class of critical Galton-Watson trees we have instead the spin glass models considered in \cite{BG} and \cite{Gun}, and the trap models with slowly-varying tails we introduce in the next section.

\subsection{One-dimensional directed trap model with slowly-varying tails}
\label{onedsec}

In this section, we describe the one-dimensional trap model with which we want to compare to our main model, and the results we are able to prove for it. To start with a formal definition, let $\tau=(\tau_x)_{x\in\mathbb{Z}}$ be a family of independent and identically-distributed strictly positive (and finite) random variables whose distribution has a slowly-varying tail, in the sense described by (\ref{svartail}), built on a probability space with measure $\mathbf{P}$; the sequence $\tau=(\tau_x)_{x\in\mathbb{Z}}$ will represent the trap environment. For a fixed bias parameter $\beta>1$, the directed trap model is then the continuous-time Markov process $X=(X_t)_{t\geq 0}$ with state space $\mathbb{Z}$, given by $X_0=0$ and with jump rates
\[c(x,y):=\left\{ \begin{array}{ll}
                    \left(\frac{\beta}{\beta+1}\right)\tau_x^{-1}, & \mbox{if $y=x+1$}, \\
                    \left(\frac{1}{\beta+1}\right)\tau_x^{-1}, & \mbox{if $y=x-1$,}
                  \end{array}\right.\]
and $c(x,y)=0$ otherwise. To be more explicit, for a particular realisation of $\tau$ we will write $P^{\tau}_x$ for the law of the Markov chain with the above transition rates, started from $x$; similarly to describing $P_x^{\mathcal{T}^*}$ in the previous section, we call this the quenched law for the directed trap model. The corresponding annealed law $\mathbb{P}_x$ is obtained by integrating out the environment similarly to (\ref{annealedlaw}), i.e.
\[\mathbb{P}_x(\cdot):=\int P^\tau_x(\cdot){\rm d}\mathbf{P}.\]

In studying the rate of escape of the above directed trap model, it is our initial aim to determine the rate of growth of
\[\Delta_n:=\inf\{t\geq 0:X_t=n\},\]
that is, the hitting times of level $n$ by the process $X$. The following theorem contains our main conclusion in this direction. As in the statement at (\ref{mlim}), we define $L(x)=1/\bar{F}(x)$.

{\thm\label{onedhit} As $n\rightarrow\infty$, the laws of the processes
\[\left(\frac{1}{n}L\left(\Delta_{nt}\right)\right)_{t\geq 0}\]
under $\mathbb{P}_0$ converge weakly with respect to the Skorohod $J_1$ topology on $D([0,\infty),\mathbb{R})$ to the law of the extremal process $(m(t))_{t\geq 0}$.}
\bigskip

Similarly to \cite[Remark 2.4]{Kasa}, we note that the proof of the above result may be significantly simplified in the case when $\bar{F}$ decays logarithmically. The reason for this is that, in the logarithmic case, the hitting time $\Delta_{n}$ is very well-approximated by the maximum holding time within the first $n$ vertices, and so the functional scaling limit for $(\Delta_n)_{n\geq 0}$ can be readily obtained from a simple study of the maximum holding time process. For general slowly varying functions, the same approximation does not provide tight enough control on $\Delta_{n}$ to apply this argument, and so a more sophisticated approach is required.

As a simple corollary of Theorem \ref{onedhit}, it is also possible to obtain a scaling result for the process $X$ itself. The definition of $m^{-1}$ should be recalled from (\ref{minv}). We similarly define the right-continuous inverse $\bar{F}^{-1}$ of $\bar{F}$, only with $>$ replaced by $<$.

{\cor\label{oned} As $n\rightarrow\infty$, the laws of the processes
\[\left(\frac{1}{n}X_{\bar{F}^{-1}(1/nt)}\right)_{t\geq 0}\]
under $\mathbb{P}_0$ converge weakly with respect to the Skorohod $M_1$ topology on $D([0,\infty),\mathbb{R})$ to the law of $(m^{-1}(t))_{t\geq 0}$.}

{\rem \label{contrem} (i) Although the preceding corollary does look somewhat awkward, it becomes much clearer for concrete choices of $\bar{F}$. For example, if $\bar{F}$ has the form described at (\ref{taileg}), then the above result concerns the distributional limit of
\[\left(\frac{1}{n}X_{e^{(nt)^{1/\gamma}}}\right)_{t\geq 0}.\]
Moreover, it can be deduced from the above result that, as $t\rightarrow\infty$, the random variable $\bar{F}(t)X_t$ converges in distribution under $\mathbb{P}_0$ to $m^{-1}(1)$, which is easily checked to have a mean one exponential distribution.\\
(ii) In a number of places in the proofs of Theorem \ref{onedhit} and Corollary \ref{oned}, we are slightly cavalier about assuming that $\bar{F}(\bar{F}^{-1}(x))=x$ for $x\in (0,1)$. This is, of course, only true in general when $\bar{F}$ is continuous. In the case when this condition is not satisfied, however, we can easily overcome the difficulties that arise by replacing $\bar{F}$ with any non-increasing continuous function $\bar{G}$ that satisfies $\bar{G}(0)=1$ and $\bar{G}(u)\sim\bar{F}(u)$ as $u\rightarrow\infty$. For example, one could define such a $\bar{G}$ by setting $\bar{G}(u):=(\frac{1}{u}\int_0^u L(v)dv)^{-1}$.}
\bigskip

The extremal aging result we are able to prove in this setting is as follows.

{\thm \label{extaging} For any $0<a<b$, we have
\[\lim_{n\rightarrow\infty}\mathbb{P}_0\left(X_{\bar{F}^{-1}(1/na)}=X_{\bar{F}^{-1}(1/nb)}\right)=\frac{a}{b}.\]}

{\rem\label{extrem} Note that if $\bar{F}_n$ and $\bar{F}$ are not continuous and eventually strictly decreasing, a minor modification to the proof of the above result (cf. Remark \ref{contrem}(ii)) is needed.}

\subsection{Article outline and notes}

The remainder of the article is organised as follows. In Section \ref{onedd}, we study the one-dimensional trap model introduced in Section \ref{onedsec} above, proving Theorem \ref{onedhit} and Corollary \ref{oned}. In Section \ref{gwtsec}, we then adapt the relevant techniques to derive Theorem \ref{dngrowth} and Corollary \ref{treecor} for the Galton-Watson tree model. The arguments of both these sections depend on the extension of the limit at (\ref{mlim}) that is proved in Section \ref{iidsec}. Before this, in Section \ref{agingsec}, we derive the extremal aging results of Theorems \ref{extagingtree} and \ref{extaging}. Finally, as noted earlier, the appendix recalls some basic facts concerning Skorohod space.

We finish the introduction with some notes about the conventions used in this article. Firstly, there are two widely used versions of the geometric distribution with a given parameter, one with support $0,1,2,\dots$ and one with support $1,2,3,\dots$. In the course of this work, we will use both, and hope that, even without explanation, it is clear from the context which version applies when. Secondly, there are many instances when for brevity we use a continuous variable where a discrete argument is required, in such places $x$, say, should be read as $\lfloor x\rfloor$. Finally, we recall that $f\sim g$ will mean
$f(x)/g(x)\rightarrow 1$ as $x\rightarrow\infty$.

\section{Directed trap model with slowly-varying tails}\label{onedd}

This section is devoted to the proof of Theorem \ref{onedhit} and Corollary \ref{oned}. To this end, we start by deriving some slight adaptations of results from \cite{Zindy} regarding the trap environment. First, define a level $n$ critical depth for traps of the environment by setting
\begin{equation}\label{gndef}
g(n):=\bar{F}^{-1}(n^{-1}\ln n).
\end{equation}
We will say that there are deep traps at the sites $\mathcal{D}:=\{x\in\mathbb{Z}:\:\tau_x > g(n)\}$, and consider the following events: for $n\in\mathbb{N}$, $T\in(0,\infty)$,
\[\mathcal{E}_1(n,T):=\left\{\min_{\substack{x_1,x_2\in\mathcal{D}\cap[1,nT]:\\x_1\neq x_2}}|x_1-x_2|> n^\kappa\right\},\]
\[\mathcal{E}_2(n):=\left\{\mathcal{D}\cap [-(\ln n)^{1+\gamma},0]=\emptyset\right\},\]
where $\kappa, \gamma\in(0,1)$ are fixed. The event $\mathcal{E}_1(n,T)$ requires that the distance between any two deep traps in the interval $[1,nT]$ is large, and the event $\mathcal{E}_2(n)$ will help to ensure that the time the process $X$ spends outside of the strictly positive integers is negligible.

{\lem \label{prelim1} Fix $T\in(0,\infty)$. As $n\rightarrow \infty$, the $\mathbf{P}$-probability of the events $\mathcal{E}_1(n,T)$ and $\mathcal{E}_2(n)$ converge to one.}
\begin{proof} To check the result for $\mathcal{E}_1(n,T)$, we simply observe that
\[\mathbf{P}\left(\mathcal{E}_1(n,T)^c\right)\leq \sum_{\substack{\{x_1,x_2\}\subseteq [1,nT]:\\0<|x_1-x_2|\leq n^\kappa}}\mathbf{P}\left(\tau_{x_1},\tau_{x_2}> g(n)\right)\leq Tn^{1+\kappa}\bar{F}\left(g(n)\right)^2\leq\frac{T(\ln n)^2}{n^{1-\kappa}}
\rightarrow 0.\]
Similarly, we have that $\mathbf{P}(\mathcal{E}_2(n)^c)\leq n^{-1}(1+(\ln n)^{1+\gamma})\ln n$, which also converges to 0.
\end{proof}

We continue by introducing the embedded discrete-time random walk associated with $X$ and some of its properties, which will be useful throughout the remainder of the section. In particular, first let $S(0)=0$ and $S(n)$ be the time of the $n$th jump of $X$; this is the clock process corresponding to $X$. The embedded discrete-time random walk is then the process $Y=(Y_n)_{n\geq 0}$ defined by setting $Y_n:=X_{S(n)}$. Clearly $Y$ is a biased random walk on $\mathbb{Z}$ under $P^\tau_0$ for $\mathbf{P}$-a.e. realisation of $\tau$, and thus satisfies, ${P}^\tau_0$-a.s.,
\[\frac{Y_n}{n}\rightarrow \frac{\beta-1}{\beta+1}>0.\]
Whilst this result already tells us that the embedded random walk $Y$ drifts off to $+\infty$ and that the time it takes to hit level $n$, that is,
\[\Delta^Y_n:=\inf\{k\geq 0: Y_k=n\},\]
is finite for each $n$, $P^\tau_0$-a.s., we further require that it does not backtrack too much, in the sense that, for each $T\in(0,\infty)$,
\[\mathcal{E}_3(n,T):=\left\{\min_{0\leq i<j\leq \Delta^Y_{nT}}(Y_j-Y_i)>-(\ln n)^{1+\gamma}\right\}\]
occurs with high probability. This is the content of the following lemma, which is essentially contained in \cite[Lemma 3]{Zindy}.

{\lem \label{prelim2} Fix $T\in(0,\infty)$. As $n\rightarrow \infty$, the $\mathbb{P}_0$-probability of the event $\mathcal{E}_3(n,T)$ converges to one.}
\bigskip

Let us now introduce the total time the biased random walk $X$ spends at a site $x\in\mathbb{Z}$,
\[T_x:=\int_0^\infty\mathbf{1}_{\{X_t=x\}}dt.\]
To study this, first observe that the clock process $S=S(n)_{n\geq0}$ can be written
\[S(n)=\sum_{i=0}^{n-1}\tau_{Y_i}\mathbf{e}_i,\]
where $(\mathbf{e}_i)_{i\geq 0}$ is an independent sequence of mean one exponential random variables under $P^\tau_0$, independent of $Y$. Moreover, for $x\in\mathbb{Z}$, let $G(x)=\#\{n\geq 0:Y_n=x\}$ be the total number of visits of the embedded random walk $Y$ to $x$. By applying the fact that $Y$ is a random walk with a strictly positive bias, we have that if $x\geq 0$, then $G(x)$ has the geometric distribution with parameter $p = (\beta-1)/(\beta+1)$ (again for $\mathbf{P}$-a.e. realisation of $\tau$). It follows that $T_x$ is equal in distribution under $\mathbb{P}_0$ to the random variable
\begin{equation}\label{txexp}
\tau_x\sum_{i=1}^{G(x)}\mathbf{e}_i,
\end{equation}
which is almost-surely finite. We will use this characterisation of the distribution of $T_x$ to check that the time spent by $X$ in traps that are not deep is asymptotically negligible, in the sense described by the following event: for $n\in\mathbb{N}$, $T\in(0,\infty)$,
\[\mathcal{E}_4(n,T):=\left\{\sum_{i=0}^{\Delta_{nT}^Y-1}\tau_{Y_i}\mathbf{e}_i\mathbf{1}_{\{\tau_{Y_i}\leq g(n)\}}<\bar{F}^{-1}(n^{-1}(\ln n)^{1/2})\right\},\]
In particular, by similar arguments to \cite[Lemma 4]{Zindy}, we deduce the following.

{\lem \label{prelim3} Fix $T\in(0,\infty)$. As $n\rightarrow \infty$, the $\mathbb{P}_0$-probability of the event $\mathcal{E}_4(n,T)$ converges to one.}
\begin{proof} We start by checking that
\begin{equation}\label{order}
\mathbf{E}(\tau_0\mathbf{1}_{\{ \tau_0\leq g(n)\}} )= o(n^{-1}\bar{F}^{-1}(n^{-1}(\ln n)^{1/2})).
\end{equation}
To this end, let $\rho,\varepsilon\in(0,1)$, and observe that
\begin{eqnarray*}
\mathbf{E}(\tau_0\mathbf{1}_{\{\tau_0\leq g(n)\}} )&\leq & g(n)\sum_{j=0}^{\infty}\rho^j\mathbf{P}(\tau_0>\rho^{j+1}g(n))\\
&\leq & c_1g(n)\sum_{j=0}^{\infty}\rho^{j-\varepsilon (j+1)}\bar{F}(g(n))\\
&\leq & \frac{c_2g(n)\ln n}{n},
\end{eqnarray*}
where the second inequality is an application of the representation theorem for slowly varying functions (\cite[Theorem 1.2]{sen}, for example), which implies that, for any $\varepsilon>0$, there exists a constant $c_3\in(0,\infty)$ such that
\begin{equation}\label{polyvar}
\frac{\bar{F}(v)}{\bar{F}(u)}\leq
c_3\left(\frac{u}{v}\right)^\varepsilon,
\end{equation}
for all $0<v\leq u$. Again applying (\ref{polyvar}), we have that $g(n)\leq c_4 \bar{F}^{-1}(n^{-1}(\ln n)^{1/2}) (\ln n)^{-1/2\varepsilon}$. Hence, if $\varepsilon$ is chosen small enough, then (\ref{order}) holds as desired.

To proceed, note that, on $\mathcal{E}_3(n,T)$, we have that
\[\sum_{i=0}^{\Delta_{nT}^Y-1}\tau_{Y_i}\mathbf{e}_i\mathbf{1}_{\{\tau_{Y_i}\leq g(n)\}}
\leq \sum_{x=-(\ln n)^{1+\gamma}}^{nT-1}T_x\mathbf{1}_{\{\tau_{x}\leq g(n)\}}.\]
Consequently, because $E^\tau_0T_x=\tau_xE^\tau_0 G(x) \leq \frac{\beta+1}{\beta-1}\tau_x$, it follows that
\[E^\tau_0\left(\sum_{i=0}^{\Delta_{nT}^Y-1}\tau_{Y_i}\mathbf{e}_i\mathbf{1}_{\{\tau_{Y_i}\leq g(n)\}}\mathbf{1}_{\mathcal{E}_3(n,T)}\right)\leq \frac{\beta+1}{\beta-1}\sum_{x=-(\ln n)^{1+\gamma}}^{nT-1}\tau_x\mathbf{1}_{\{\tau_x\leq g(n)\}}.\]
Combining this bound with (\ref{order}) and using Markov's inequality yields
\[\mathbb{P}_0\left(\mathcal{E}_3(n,T)\cap\mathcal{E}_4(n,T)^c\right)
\leq\frac{\beta+1}{(\beta-1)\bar{F}^{-1}(n^{-1}(\ln n)^{1/2})}\mathbf{E}\left(\sum_{x=-(\ln n)^{1+\gamma}}^{nT-1}\tau_i\mathbf{1}_{\{\tau_x\leq g(n)\}}\right)=o(1).\]
On recalling the conclusion of Lemma \ref{prelim2}, this completes the proof.
\end{proof}

As a consequence of the previous result, to deduce a scaling limit for the sequence $(\Delta_n)_{n\geq0}$, it will suffice to study sums of the form $\sum_{x=1}^{n}T_x\mathbf{1}_{\{\tau_x> g(n)\}}$. In fact, the backtracking result of Lemma \ref{prelim2} will further allow us to replace $T_x$ in this expression by
\begin{equation}\label{tildetxdef}
\tilde{T}_x:=\int_{\Delta_x}^{\Delta_{x,(\ln x)^{1+\gamma}}}\mathbf{1}_{\{X_t=x\}}dt,
\end{equation}
where $\Delta_{x,(\ln x)^{1+\gamma}}$ is the first time after $\Delta_x$ that $X$ leaves the interval $[x-(\ln x)^{1+\gamma},x+(\ln x)^{1+\gamma}]$. This is particularly useful because, by applying the fact that deep traps are separated by a distance that is polynomial in $n$ (see Lemma \ref{prelim1}), it will be possible to decouple the random variables $(\tilde{T}_x\mathbf{1}_{\{\tau_x> g(n)\}})_{x\geq1}$ in such a way that enables us to deduce functional scaling results for their sums from those for independent sums proved in Section \ref{iidsec}. Before commencing this program in Lemma \ref{tildetxsums}, however, we derive a preliminary lemma that suitably describes the asymptotic behaviour of the distributional tail
\[\bar{F}_n(u):=\mathbb{P}_0\left(\tilde{T}_x\mathbf{1}_{\{\tau_x> g(n)\}}>u\right).\]
(Clearly, the definition of $\bar{F}_n$ is independent of the particular $x\geq 1$ considered.)

{\lem \label{taill} For every $\varepsilon>0$, there exists a constant $c$ such that, for any $u\geq c(g(n)\vee 1)$,
\[(1-\varepsilon)\bar{F}_n(u)\leq \bar{F}(u)\leq (1+\varepsilon)\bar{F}_n(u).\]}
\begin{proof} For $x\geq 1$, let $\tilde{G}(x)$ be the total number of visits of the embedded random walk $Y$ to $x$ up until the first time after $\Delta_x^Y$ that it leaves the interval $[x-(\ln n)^{1+\gamma},x+(\ln n)^{1+\gamma}]$. Then, similarly to (\ref{txexp}), we have that $\tilde{T}_x$ is distributed as $\tau_x\sum_{i=1}^{\tilde{G}(x)}\mathbf{e}_i$. Hence, setting $\Gamma:=\sum_{i=1}^{\tilde{G}(x)}\mathbf{e}_i$, we can use the independence of $\Gamma$ and $\tau_x$ under $\mathbb{P}_0$ to write
\begin{eqnarray*}
\bar{F}_n(u) &=&\mathbb{P}_0\left(\tau_x\Gamma >u,\:\tau_x>g(n)\right)\\
&=&\int_{0}^{u/g(n)}\bar{F}\left(uv^{-1}\right)\mathbb{P}_0\left(\Gamma\in dv\right)+\int_{u/g(n)}^\infty\bar{F}\left(g(n)\right)\mathbb{P}_0\left(\Gamma\in dv\right).
\end{eqnarray*}
It follows that
\begin{equation}\label{ud}
\left|\frac{\bar{F}_n(u)}{\bar{F}(u)}-1\right|\leq
\int_{0}^{\infty}\left|\frac{\bar{F}\left(uv^{-1}\right)}{\bar{F}(u)}-1\right|\mathbb{P}_0\left(\Gamma\in dv\right)+\left(\frac{\bar{F}\left(g(n)\right)}{\bar{F}(u)}+1\right)\mathbb{P}_0\left(\Gamma \geq u/g(n)\right).
\end{equation}

The first term on the right-hand side of (\ref{ud}) is independent of $n$, and so it will be enough for our purposes to show that it converges to 0 as $u\rightarrow\infty$. To do this, first note that, by the monotonicity of $\bar{F}$, \eqref{svartail} holds uniformly for $v\in[v_0,v_1]$ for any $0<v_0\leq v_1<\infty$. Hence, the $\limsup$ as $u\rightarrow\infty$ of the term of interest is bounded above by
\[\mathbb{P}_0\left(\Gamma\not\in[v_0,v_1]\right)+\limsup_{u\rightarrow\infty}\int_{0}^{v_0} \frac{\bar{F}(uv^{-1})}{\bar{F}(u)}\mathbb{P}_0\left(\Gamma\in dv\right)+\limsup_{u\rightarrow\infty}\int_{v_1}^\infty \frac{\bar{F}(uv^{-1})}{\bar{F}(u)}\mathbb{P}_0\left(\Gamma\in dz\right),\]
for any $0<v_0\leq v_1<\infty$. Now, if $v_0<1$, then $\bar{F}(uv^{-1})\leq \bar{F}(u)$ for all $v\in [0,v_0]$, and so the first limsup is bounded above by $\mathbb{P}_0(\Gamma\leq v_0)$. Furthermore, if $v_1$ is chosen to be no less than 1, then we can apply the bound at (\ref{polyvar}) to estimate $\bar{F}(uv^{-1})/\bar{F}(u)$ by $cv^\varepsilon$ for $v\geq v_1$. Thus
\[\limsup_{u\rightarrow\infty}
\int_{0}^{\infty}\left|\frac{\bar{F}\left(uv^{-1}\right)}{\bar{F}(u)}-1\right|\mathbb{P}_0\left(\Gamma\in dv\right)\leq 2\mathbb{P}_0\left(\Gamma\not\in[v_0,v_1]\right)+c\int_{v_1}^\infty v^\varepsilon \mathbb{P}_0\left(\Gamma\in dv\right).\]
Since $\mathbb{E}_0(\Gamma^\varepsilon)\leq \mathbb{E}_0(1+\Gamma)= 1+\mathbb{E}_0(\tilde{G}(x))\mathbb{E}_0(\mathbf{e}_1)<\infty$, by taking $v_0$ arbitrarily small and $v_1$ arbitrarily large, the upper bound here can be made arbitrarily small, meaning that
\[\lim_{u\rightarrow\infty}
\int_{0}^{\infty}\left|\frac{\bar{F}\left(uv^{-1}\right)}{\bar{F}(u)}-1\right|\mathbb{P}_0\left(\Gamma\in dv\right)=0,\]
as desired.

For the second term on the right-hand side of (\ref{ud}), we apply (\ref{polyvar}) and Markov's inequality to deduce that, if $u\geq g(n)$, then
\[\left(\frac{\bar{F}\left(g(n)\right)}{\bar{F}(u)}+1\right)\mathbb{P}_0\left(\Gamma \geq u/g(n)\right)\leq \left(c \left(\frac{u}{g(n)}\right)^\varepsilon+1\right)\frac{g(n)\mathbb{E}_0\Gamma}{u},\]
where $\varepsilon\in(0,1)$ is fixed. Thus, since the above bound is small whenever $u/g(n)$ is large (it was already noted in the previous paragraph that $\Gamma$ has a finite first moment), the proof is complete.
\end{proof}

{\lem \label{tildetxsums} As $n\rightarrow\infty$, the laws of the processes
\[\left(\frac{1}{n}L\left(\sum_{x=1}^{nt}\tilde{T}_x\mathbf{1}_{\{\tau_x> g(n)\}}\right)\right)_{t\geq0}\]
under $\mathbb{P}_0$ converge weakly with respect to the Skorohod $J_1$ topology on $D([0,\infty),\mathbb{R})$ to the law of $(m(t))_{t\geq 0}$.}
\begin{proof} First, fix $T\in(0,\infty)$ and suppose $(f_x)_{x\geq 1}$ is a collection of bounded, continuous functions on $\mathbb{R}$. We then have that
\begin{eqnarray*}
\lefteqn{\mathbb{E}_0\left(\mathbf{1}_{\mathcal{E}_1(n,T)}\prod_{x=1}^{nT} f_x\left(\tilde{T}_x\mathbf{1}_{\{\tau_x> g(n)\}}\right)\right)}\\
&=&\sum_{B}\mathbb{E}_0\left(\mathbf{1}_{\{\mathcal{D}\cap[1,nT]=B\}}
\prod_{x=1}^{nT} f_x\left(\tilde{T}_x\mathbf{1}_{\{\tau_x> g(n)\}}\right)\right)\\
&=&\sum_B\mathbb{E}_0\left(\prod_{x\in B}f_x\left(\tilde{T}_x\mathbf{1}_{\{\tau_x>g(n)\}}\right)\mathbf{1}_{\{\tau_x>g(n)\}}\prod_{x\in [1,nT]\backslash B}f_x(0)\mathbf{1}_{\{\tau_x\leq g(n)\}}\right),
\end{eqnarray*}
where the sums are over subsets $B\subseteq [1,nT]$ such that if $x_1,x_2\in B$ and $x_1\neq x_2$, then $|x_1-x_2|> n^\kappa$. By applying the independence of traps at different sites and the disjointness of the intervals $([x-(\ln n)^{1+\gamma},x+(\ln n)^{1+\gamma}])_{x\in B}$ for the relevant choices of $B$, the above sum can be rewritten as
\[\sum_B\prod_{x\in B}\mathbb{E}_0\left(f_x\left(\tilde{T}_x\mathbf{1}_{\{\tau_x>g(n)\}}\mathbf{1}_{\{\tau_x>g(n)\}}\right)\right)
\prod_{x\in [1,nT]\backslash B}\mathbb{E}_0\left(f_x(0)\mathbf{1}_{\{\tau_x\leq g(n)\}}\right).\]
In particular, it follows that
\[\mathbb{E}_0\left(\mathbf{1}_{\mathcal{E}_1(n,T)}\prod_{x=1}^{nT} f_x\left(\tilde{T}_x\mathbf{1}_{\{\tau_x> g(n)\}}\right)\right)
=\mathbb{E}_0\left(\mathbf{1}_{\mathcal{E}'_1(n,T)}\prod_{x=1}^{nT} f_x\left(\tilde{T}'_x\mathbf{1}_{\{\tau'_x> g(n)\}}\right)\right),\]
where we suppose that, under $\mathbb{P}_0$, the pairs of random variables $(\tilde{T}'_x,\tau'_x)$, ${x\geq1}$, are independent and identically-distributed as $(\tilde{T}_1,\tau_1)$, and the event $\mathcal{E}'_1(n,T)$ is defined analogously to $\mathcal{E}_1(n,T)$ from these random variables. Consequently, under $\mathbb{P}_0$, the laws of $(\tilde{T}_x\mathbf{1}_{\{\tau_x> g(n)\}})_{x=1}^{nT}$ conditional on $\mathcal{E}_1(n,T)$ and
$(\tilde{T}'_x\mathbf{1}_{\{\tau'_x> g(n)\}})_{x=1}^{nT}$ conditional on $\mathcal{E}'_1(n,T)$ are identical.

By applying the conclusion of the previous paragraph, we obtain that, for any bounded function $H:D([0,T],\mathbb{R})\rightarrow \mathbb{R}$ that is continuous with respect to the Skorohod $J_1$ topology,
\[\left|\mathbb{E}_0\left[H\left(\left(\frac{1}{n}L\left(\sum_{x=1}^{nt}\tilde{T}_x\mathbf{1}_{\{\tau_x> g(n)\}}\right)\right)_{t\in[0,T]}\right)\right]-\mathbb{E}_0\left[H\left(\left(\frac{1}{n}L\left(\sum_{x=1}^{nt}\tilde{T}'_x\mathbf{1}_{\{\tau'_x> g(n)\}}\right)\right)_{t\in[0,T]}\right)\right]\right|\]
is bounded above by $2\|H\|_\infty\mathbf{P}\left(\mathcal{E}_1(n,T)^c\right)$. Since Lemma \ref{prelim1} tells us that this upper bound converges to 0 as $n\rightarrow\infty$, to complete the proof it will thus suffice to establish the result with $(\tilde{T}_x,\tau_x)_{x\geq1}$ replaced by $(\tilde{T}'_x,\tau'_x)_{x\geq1}$. However, because we are assuming that $(\tilde{T}'_x,\tau'_x)_{x\geq1}$ are independent, the tail asymptotics proved in Lemma \ref{taill} allow us to derive the relevant scaling limit for the sums involving $(\tilde{T}'_x,\tau'_x)_{x\geq1}$ by a simple application of Theorem \ref{thm:Kasa-ext} (with $h_1(n)=\ln n$ and $h_2(n)=0$).
\end{proof}

We are now in a position to prove Theorem \ref{onedhit} by showing that the rescaled sums considered in the previous lemma suitably well approximate the sequence $(\Delta_n)_{n\geq 1}$.

\begin{proof}[Proof of Theorem \ref{onedhit}] Fix $T\in(0,\infty)$ and observe that, on $\mathcal{E}_2(n)\cap\mathcal{E}_3(n,T)\cap\mathcal{E}_4(n,T)$, we have that
\begin{equation}\label{dnbound}
\sum_{x=1}^{nt-(\ln n)^{1+\gamma}}\tilde{T}_x\mathbf{1}_{\{\tau_x>g(n)\}}\leq \Delta_{nt}\leq \sum_{x=1}^{nt}\tilde{T}_x\mathbf{1}_{\{\tau_x>g(n)\}}+\bar{F}^{-1}(n^{-1}(\ln n)^{1/2}),\hspace{20pt}\forall t\in[0,T].
\end{equation}
By reparameterising the time-scales in the obvious way, it is clear that
\begin{equation}\label{dj1bound}
d_{J_1}\left(\left(\frac{1}{n}L\left(\sum_{x=1}^{nt-(\ln n)^{1+\gamma}}\tilde{T}_x\mathbf{1}_{\{\tau_x>g(n)\}}\right)\right)_{t\in[0,T]},
\left(\frac{1}{n}L\left(\sum_{x=1}^{nt}\tilde{T}_x\mathbf{1}_{\{\tau_x>g(n)\}}\right)\right)_{t\in[0,T]}\right),
\end{equation}
where $d_{J_1}$ is the Skorohod $J_1$ distance on $D([0,T],\mathbb{R})$ (as defined in the appendix at (\ref{dj1})), is bounded above by
\[\frac{(\ln n)^{1+\gamma}}{n}+\frac{1}{n}L\left(\sum_{x=1}^{nT}\tilde{T}_x\mathbf{1}_{\{\tau_x>g(n)\}}\right)-\frac{1}{n}L\left(\sum_{x=1}^{n(T-\varepsilon)}\tilde{T}_x\mathbf{1}_{\{\tau_x>g(n)\}}\right),\]
for large $n$. (Note that the first term above relates to the distortion of the time scale needed to compare the two processes.) By Lemma \ref{tildetxsums}, this bound converges in distribution under $\mathbb{P}_0$ to $m(T)-m(T-\varepsilon)$. Now, in the limit as $\varepsilon\rightarrow 0$,
$m(T)-m(T-\varepsilon)$ converges to 0 in probability. It readily follows that, as $n\rightarrow\infty$, so does the expression at (\ref{dj1bound}). Hence, the theorem will follow from Lemmas \ref{prelim1}, \ref{prelim2}, \ref{prelim3} and \ref{tildetxsums}, if we can show that
\[\sup_{x\in[0,\Xi]}\frac{1}{n}\left|L\left(x+\bar{F}^{-1}(n^{-1}(\ln n)^{1/2})\right)-L(x)\right|\rightarrow 0,\]
in $\mathbb{P}_0$ probability, where $\Xi:=\sum_{x=1}^{nT}\tilde{T}_x\mathbf{1}_{\{\tau_x>g(n)\}}$. To check this, we start by noting that Lemma \ref{tildetxsums} implies, for any $\lambda>0$,
\[\mathbb{P}_0\left(\Xi\leq \bar{F}^{-1}(1/n\lambda)\right)= \mathbb{P}_0\left(n^{-1}L(\Xi)\leq \lambda\right)\rightarrow \mathbb{P}_0(m(T)\leq \lambda).\]
By choosing $\lambda$ suitably large, the limiting probability can be made arbitrarily close to 1. Thus the problem reduces to showing that, for any $\lambda\in(0,\infty)$,
\[\sup_{x\in[0,\bar{F}^{-1}(1/n\lambda)]}\frac{1}{n}\left|L\left(x+\bar{F}^{-1}(n^{-1}(\ln n)^{1/2})\right)-L(x)\right|\rightarrow 0.\]
Let $\varepsilon\in(0,\lambda)$, then, since $\bar{F}^{-1}(n^{-1}(\ln n)^{1/2})\leq \bar{F}^{-1}(1/n\varepsilon)$ for large enough $n$, we have that
\begin{eqnarray*}
\lefteqn{\sup_{x\in[\bar{F}^{-1}(1/n\varepsilon),\bar{F}^{-1}(1/n\lambda)]}\frac{1}{n}\left|L\left(x+\bar{F}^{-1}(n^{-1}(\ln n)^{1/2})\right)-L(x)\right|}\\
&\leq &\frac{1}{n}L\left(\bar{F}^{-1}(1/n\lambda)\right)\sup_{x\geq\bar{F}^{-1}(1/n\varepsilon)}\left|\frac{L(x+\bar{F}^{-1}(n^{-1}(\ln n)^{1/2}))}{L(x)}-1\right|\\
&\leq &\lambda\sup_{x\geq\bar{F}^{-1}(1/n\varepsilon)}\left|\frac{L(2x)}{L(x)}-1\right|,
\end{eqnarray*}
which converges to 0 as $n\rightarrow\infty$ by (\ref{svartail}). Moreover, we also have that
\[\sup_{x\in[0,\bar{F}^{-1}(1/n\varepsilon)]}\frac{1}{n}\left|L\left(x+\bar{F}^{-1}(n^{-1}(\ln n)^{1/2})\right)-L(x)\right|
\leq  \frac{1}{n}L(2\bar{F}^{-1}(1/n\varepsilon))
\sim  \frac{1}{n}L(\bar{F}^{-1}(1/n\varepsilon)),
\]
where the asymptotic equivalence is an application of (\ref{svartail}). In particular, since the right-hand side above is equal to $\varepsilon$, which can be  chosen arbitrarily small, the result follows.
\end{proof}

From this, the proof of Corollary \ref{oned} is relatively straightforward.

\begin{proof}[Proof of Corollary \ref{oned}] Define $X^*=(X_t)_{t\geq0}$ to be the running supremum of $X$, i.e. $X_t^*:=\max_{s\leq t}X_s$. Since $X_t^*\geq n$ if and only if $\Delta_n\leq t$, we obtain that $(X^*_t+1)_{t\geq 0}$ is the inverse of $(\Delta_n)_{n\geq 0}$ in the sense described at (\ref{minv}). Thus, because the inverse map is continuous with respect to the Skorohod $M_1$ topology (at least on the subset of functions $f\in D([0,\infty),\mathbb{R})$ that satisfy $\limsup_{t\rightarrow\infty} f(t)=\infty$, see \cite{Whitt}), it is immediate from Theorem \ref{onedhit} that, as $n\rightarrow\infty$, the laws of the processes
\[\left(\frac{1}{n}X^*_{\bar{F}^{-1}(1/nt)}\right)_{t\geq 0}\]
under $\mathbb{P}_0$ converge weakly with respect to the Skorohod $M_1$ topology on $D([0,\infty),\mathbb{R})$ to the law of $(m^{-1}(t))_{t\geq 0}$. Thus, to complete the proof, it will suffice to demonstrate that, for any $T\in(0,\infty)$,
\[\sup_{t\in[0,T]}\frac{1}{n}\left|X^*_{\bar{F}^{-1}(1/nt)}-X_{\bar{F}^{-1}(1/nt)}\right|\rightarrow0\]
in $\mathbb{P}_0$-probability as $n\rightarrow\infty$. To do this, we first fix $T\in(0,\infty)$ and set $N:=nT\ln (nT)$. Theorem \ref{onedhit} then implies that $\mathbb{P}_0\left(\Delta_{N}\geq \bar{F}^{-1}(1/nT)\right)\rightarrow 1$ as $n\rightarrow \infty$. Moreover, on the set $\{\Delta_{N}\geq \bar{F}^{-1}(1/nT)\}$, it is the case that
\[\sup_{t\in[0,T]}\left|X^*_{\bar{F}^{-1}(1/nt)}-X_{\bar{F}^{-1}(1/nt)}\right|\leq\sup_{k\leq \Delta^Y_{N} }(Y^*_k-Y_k),\]
where $Y^*$ is the running supremum of $Y$. Hence
\begin{eqnarray*}
\lefteqn{\limsup_{n\rightarrow\infty}\mathbb{P}_0\left(\sup_{t\in[0,T]}\frac{1}{n}\left|X^*_{\bar{F}^{-1}(1/nt)}-X_{\bar{F}^{-1}(1/nt)}\right|>\varepsilon\right)}\\
&\leq&
\limsup_{n\rightarrow\infty}\mathbb{P}_0\left(\frac{1}{n}\sup_{k\leq \Delta^Y_{N} }(Y^*_k-Y_k)>\varepsilon\right)\\
&\leq&\limsup_{n\rightarrow\infty}\mathbb{P}_0\left(n^{-1}(\ln N)^{1+\gamma}>\varepsilon,\mathcal{E}_3(N,1)\right)\\
&=&0,
\end{eqnarray*}
where we have applied the fact that $\mathbb{P}_0(\mathcal{E}_3(N,1)^c)\rightarrow0$, which is the conclusion of Lemma \ref{prelim2}, and also that $n^{-1}(\ln N)^{1+\gamma}\rightarrow 0$, which is clear from the definition of $N$.
\end{proof}

\section{Biased random walk on critical Galton-Watson trees}\label{gwtsec}

In this section, we explain how techniques similar to those of the previous section can be used to deduce the corresponding asymptotics for a biased random walk on a critical Galton-Watson tree conditioned to survive. Prior to proving our main results (Theorem \ref{dngrowth} and Corollary \ref{treecor}), however, we proceed in the next two subsections to derive certain properties regarding the structure of the tree $\mathcal{T}^*$ and deduce some preliminary simple random walk estimates, respectively. These results establish information in the present setting that is broadly analogous to that contained in Lemmas \ref{prelim1}-\ref{taill} for the directed trap model.

\subsection{Structure of the infinite tree}\label{strucsec}

\begin{figure}[ht]
\begin{center}
\scalebox{0.35}{\includegraphics{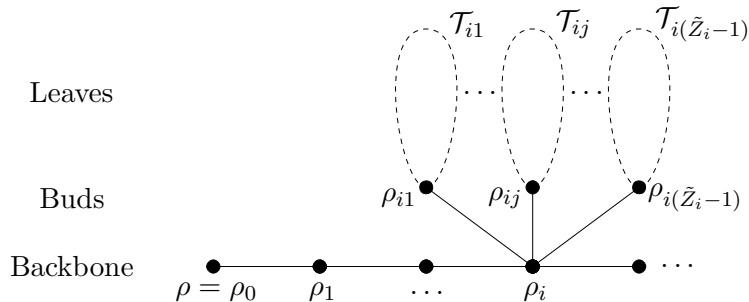}}
\rput(-6.1,-.25){$\rho=\rho_0$}
\rput(-4.7,-.25){$\rho_1$}
\rput(-3.3,-.25){\dots}
\rput(-1.9,-.25){$\rho_i$}
\rput(0,0.1){\dots}
\rput(-3.7,1.0){$\rho_{i1}$}
\rput(-2.3,1.0){$\rho_{ij}$}
\rput(0.2,1.0){$\rho_{i(\tilde{Z}_i-1)}$}
\rput(-2.6,2.4){\dots}
\rput(-1.2,2.4){\dots}
\rput(-8,0.1){Backbone}
\rput(-8,1.0){Buds}
\rput(-8,2.4){Leaves}
\rput(-2.8,3.3){$\mathcal{T}_{i1}$}
\rput(-1.4,3.3){$\mathcal{T}_{ij}$}
\rput(0.3,3.3){$\mathcal{T}_{i(\tilde{Z}_i-1)}$}
\end{center}
\caption{Decomposition of $\mathcal{T}^*$.}
\label{treedecomp}
\end{figure}

A key tool throughout this study is the spinal decomposition of $\mathcal{T}^*$ that appears as \cite[Lemma 2.2]{Kesten}, and which can be described as follows. First, $\mathbf{P}$-a.e. realisation of $\mathcal{T}^*$ admits a unique non-intersecting infinite path starting at the root. Conditional on this `backbone', the number of children of vertices on the backbone are independent, each distributed as a size-biased random variable $\tilde{Z}$, which satisfies
\begin{equation}\label{sizebiasz}
\mathbf{P}\left(\tilde{Z}=k\right)=k\mathbf{P}(Z=k),\hspace{20pt}k\geq 1.
\end{equation}
Moreover, conditional on the backbone and the number of children of each backbone element, the trees descending from the children of backbone vertices that are not on the backbone are independent copies of the original critical branching process $\mathcal{T}$. To fix notation and terminology for this decomposition, we will henceforth suppose that $\mathcal{T}^*$ has been built by starting with a semi-infinite path, $\{\rho=\rho_0,\rho_1,\rho_2,\dots\}$ -- this will form the backbone of $\mathcal{T}^*$. Then, after selecting $(\tilde{Z}_i)_{i\geq 0}$ independently with distribution equal to that of $\tilde{Z}$, to each backbone vertex $\rho_i$, we attach a collection of `buds' $\rho_{ij}$, $j=1,\dots,\tilde{Z}_i-1$. Finally, we grow from each bud $\rho_{ij}$ a `leaf' $\mathcal{T}_{ij}$, that is, a Galton-Watson tree with initial ancestor $\rho_{ij}$ and offspring distribution $Z$. See Figure \ref{treedecomp} for a graphical representation of these definitions.

With this picture, it is clear how we can view $\mathcal{T}^*$ as an essentially one-dimensional trap model with the backbone playing the role of $\mathbb{Z}$ in the previous section. Rather than having an exponential holding time at each vertex $\rho_i$, however, we have a random variable representing the time it takes $X$ to leave the tree $\mathcal{T}_{i}:=\{\rho_i\}\cup(\cup_{j=1,\dots,\tilde{Z}_i-1}\mathcal{T}_{ij})$ starting from $\rho_i$. As will be made precise later, key to determining whether this time is likely to be large or not are the heights of the leaves connected to $\rho_i$. For this reason, the rest of this section will be taken up with an investigation into the big, or perhaps more accurately tall, leaves of $\mathcal{T}^*$.

More concretely, we start by introducing a sequence of critical heights $(h_n)_{n\geq 1}$ by setting $h_n:=n(\ln n)^{-1}$ (roughly, $\beta^{h_n}$ will play the role that the $g(n)$ introduced at (\ref{gndef}) did in the previous section), and define, for each $i\geq0$,
\[N_n(i):=\#\left\{1\leq j\leq \tilde{Z}_i-1:\:h(\mathcal{T}_{ij})\geq h_n\right\},\]
where $h(\mathcal{T}_{ij})$ is the height of the tree $\mathcal{T}_{ij}$, so that $N_n(i)$ counts the number of big leaves emanating from the backbone vertex $\rho_i$. The random variables in the collection $(N_n(i))_{i\geq 0}$ are independent and identically-distributed. Moreover, it is possible to describe the asymptotic probability that one of these random variables is equal to zero, i.e. there is no big leaf at the relevant site.

{\lem\label{bigtrap} Let $\alpha\in(1,2]$. As $n\rightarrow\infty$, we have that
\[\mathbf{P}\left(N_n(0)=0\right)\sim 1-\frac{\alpha}{(\alpha-1)h_n}.\]}
\begin{proof} By conditioning on the number of buds attached to the root, we have
\[\mathbf{P}\left(N_n(0)=0\right)=\mathbf{E}\left(\left(1-q_{h_n}\right)^{\tilde{Z}-1}\right),\]
where, as introduced above (\ref{probdecay}), $q_k$ is the probability that an unconditioned branching process with offspring distribution $Z$ survives for at least $k$ generations. By the size-biasing of (\ref{sizebiasz}), this can be rewritten as
\[\mathbf{P}\left(N_n(0)=0\right)=\mathbf{E}\left(Z\left(1-q_{h_n}\right)^{{Z}-1}\right)=f'\left(1-q_{h_n}\right),\]
where $f'$ is the derivative of the generating function $f$, as defined at (\ref{genfun}). Now, by \cite[(2.1)]{Slack}, it holds that $f'(1-x)\sim 1-\alpha x^{\alpha-1}L(x)$ as $x\rightarrow 0^+$, and so
\[\mathbf{P}\left(N_n(0)=0\right)\sim 1- \alpha q_{h_n}^{\alpha-1}L(q_{h_n}).\]
From this, the proof is completed by recalling the tail decay at (\ref{probdecay}).
\end{proof}

It will be important for our future arguments that the sites from which big leaves emanate are not too close together, and that there are no big traps close to $\rho$. The final lemma of this section demonstrates that the sequence of critical heights we have chosen achieves this.

{\lem\label{separation} Let $\alpha\in(1,2]$, $T\in(0,\infty)$ and $\varepsilon\in(0,1)$. As $n\rightarrow\infty$,
\[\mathbf{P}\left(\sum_{i=m}^{m+n^\varepsilon}\mathbf{1}_{\{N_n(i)\geq 1\}}\geq 2\mbox{ for some }m\in\{0,1,\dots,Tn-n^\varepsilon\}\right)\rightarrow0\]
and also
\[\mathbf{P}\left(\sum_{i=0}^{n^\varepsilon}\mathbf{1}_{\{N_n(i)\geq 1\}}\geq 1\right)\rightarrow0.\]}
\begin{proof} This is essentially the same as Lemma \ref{prelim1}.
\end{proof}

\subsection{Initial random walk estimates}\label{rwsec}

This section collects together some preliminary results for the biased random walk $(X_m)_{m\geq 0}$ on $\mathcal{T}^*$, regarding in particular: the amount of backtracking performed by the embedded biased random walk on the backbone; the amount of time $X$ spends in small leaves; the amount of time $X$ spends close to the base of big leaves; and tail estimates for the amount of time $X$ spends deep within big leaves.

To begin with, we introduce $Y=(Y_n)_{n\geq 1}$ to represent the jump process of $\pi(X)$, where $\pi:\mathcal{T}^*\rightarrow \{\rho_0,\rho_1,\dots\}$ is the projection onto the backbone, i.e. $\pi(x)=\rho_i$ for $x\in \mathcal{T}_i$. More precisely, set $S(0)=0$,
\[S(n)=\inf\left\{m> S(n-1):\pi(X_{m})\neq \pi(X_{m-1})\right\},\hspace{20pt}\forall n\geq 1,\]
and then define $Y_n:=X_{S(n)}$. From this construction, it is clear that, under either the quenched or annealed law, $Y$ is simply a biased random walk on the semi-infinite line graph $\{\rho_0,\rho_1,\dots\}$, and so, as in the previous section, we can control the amount it backtracks. In particular, if we let
\[
\Delta_n^Y:=\inf\left\{m\geq 0:\:Y_m=\rho_n\right\}
\]
be the first time that the embedded random walk $Y$ reaches level $n$ along the backbone, then we have the following result, which is simply a restatement of Lemma \ref{prelim2}. We recall that $d_{\mathcal{T}^*}$ is the shortest path graph distance on $\mathcal{T}^*$.

{\lem\label{backtrack} Let $\alpha\in(1,2]$, $T\in(0,\infty)$ and $\gamma>0$. As $n\rightarrow\infty$,
\[\mathbb{P}_\rho\left(\min_{0\leq i<j\leq \Delta^Y_{nT}}\left(d_{\mathcal{T}^*}(\rho_0,Y_j)-d_{\mathcal{T}^*}(\rho_0,Y_i)\right)\leq -(\ln n)^{1+\gamma}\right)\rightarrow 0.\]}

Our next goal is to show that the time the biased random walk $X$ spends outside of the big leaves of $\mathcal{T}^*$ is unimportant, where we define the set of vertices in big leaves to be
\[\mathcal{B}:=\left\{x\in \mathcal{T}_{ij}:\:i\geq0,\: 1\leq j\leq{\tilde{Z}_i-1},\:h(\mathcal{T}_{ij})\geq h_n\right\}.\]
Key to doing this is the following equality, which is obtained by applying standard results for weighted random walks on graphs (cf. \cite[Lemma 3.1]{AFGH}):
\begin{equation}\label{treeexit}
E^{\mathcal{T}^*}_{\rho_{ij}}\tau_{\rho_i}=1+\beta^{-i}\sum_{\substack{x,y\in \mathcal{T}_{ij}:\\ x\sim y}}c(x,y),
\end{equation}
where for a vertex $x\in\mathcal{T}^*$, we define $\tau_{x}:=\inf\{m\geq 0:\:X_m=x\}$. For the statement of the next lemma, which is approximately analogous to Lemma \ref{prelim3}, we recall the definition of $(\Delta_n)_{n\geq0}$ from (\ref{xhit}).

{\lem \label{smallleaf} Let $\alpha\in (1,2]$, $T\in(0,\infty)$ and $\varepsilon>0$.  As $n\rightarrow\infty$,
\[\mathbb{P}_\rho\left(\sum_{m\leq \Delta_{nT}}\mathbf{1}_{\{X_m\not\in\mathcal{B}\}}\geq \beta^{h_n(1+\varepsilon)} \right)\rightarrow 0.\]}
\begin{proof} We start by estimating the quenched expectation of the time $X$ spends in a particular small leaf before reaching level $nT$ along the backbone. Thus, suppose we have a leaf $\mathcal{T}_{ij}$ such that $i<nT$ and $h(\mathcal{T}_{ij})<h_n$. Starting from the vertex $\rho_i$, the probability of hitting $\rho_{ij}$ before $\rho_{nT}$ can be computed exactly, by elementary means, as
\[P^{\mathcal{T}^*}_{\rho_i}\left(\tau_{\rho_{ij}}<\tau_{\rho_{nT}}\right)=\frac{1+\beta^{-1}+\dots+\beta^{i+1-nT}}{1+(1+\beta^{-1}+\dots+\beta^{i+1-nT})}\leq \frac{1}{2-\beta^{-1}}.\]
This means that the number of separate visits $X$ makes to $\mathcal{T}_{ij}$ is stochastically dominated by a geometric random variable with parameter $1-(2-\beta^{-1})^{-1}$, and so its mean is bounded above by $\beta/(\beta-1)$. Moreover, the equality at (\ref{treeexit}) and our assumption on $h(\mathcal{T}_{ij})$ imply that, on each visit to $\mathcal{T}_{ij}$, the amount of time $X$ spends there is bounded above by
\[E^{\mathcal{T}^*}_{\rho_{ij}}\tau_{\rho_i}\leq 1+2\beta^{h_n}\#\mathcal{T}_{ij},\]
where $\#\mathcal{T}_{ij}$ is the total number of vertices in $\mathcal{T}_{ij}$. Hence
\begin{equation}\label{est1}
E^{\mathcal{T}^*}_{\rho}\left(\sum_{m\leq \Delta_{nT}}\mathbf{1}_{\{X_m\in\mathcal{T}_{ij}\}}\right)\leq \frac{\beta}{\beta-1}\left(1+2\beta^{h_n}\#\mathcal{T}_{ij}\right).
\end{equation}

As for the estimating time spent at a vertex $\rho_i$, where $0<i<nT$, we start by noting that the total number of returns to $\rho_i$ is a geometric random variable. Moreover, its parameter $P_{\rho_i}^{\mathcal{T}^*}(\tau_{\rho_i}^+=\infty)$, where $\tau_{\rho_i}^+:=\inf\{m>0:X_m=\rho_i\}$ is the first return time to $\rho_i$, can easily be bounded below by the probability that $X$ jumps from $\rho_i$ to $\rho_{i+1}$ on its first step times the probability that a biased random walk on $\mathbb{Z}$ never hits the vertex to the left of its starting point. Since the first of these quantities is given by $\beta/(\beta\tilde{Z}_i+1)$ and the second is equal to $1-\beta^{-1}$, it follows that
\begin{equation}\label{est2}
E_{\rho}^{\mathcal{T}^*}\left(\sum_{m\leq \Delta_{nT}} \mathbf{1}_{\{X_m=\rho_i\}}\right)\leq c_{\beta} \tilde{Z}_i.
\end{equation}
A similar argument applies for $i=0$.

Piecing together the estimates at (\ref{est1}) and (\ref{est2}), we thus obtain
\begin{equation}\label{ebound}
E^{\mathcal{T}^*}_{\rho}\left(\sum_{m\leq \Delta_{nT}}\mathbf{1}_{\{X_m\not\in\mathcal{B}\}}\right)\leq c_\beta\beta^{h_n}\sum_{i=0}^{nT-1}\tilde{Z}_i\left[1+\max_{j=1,\dots,\tilde{Z}_i-1}\#\mathcal{T}_{ij}\right],
\end{equation}
where $c_\beta$ is a constant depending only on $\beta$. Now, to bound the summands, we consider the following probabilistic upper bound
\begin{equation}\label{twoterms}
\mathbf{P}\left(\tilde{Z}_i\left[1+\max_{j=1,\dots,\tilde{Z}_i-1}\#\mathcal{T}_{ij}\right]\geq k\right)\leq \mathbf{P}\left(\tilde{Z}_i\geq k^{1/2}\right)+\mathbf{P}\left(\max_{j=1,\dots,\tilde{Z}_i-1}\#\mathcal{T}_{ij}\geq k^{1/2}-1\right).
\end{equation}
For the first of these terms, we apply the size-biasing of (\ref{sizebiasz}) and Markov's inequality to deduce
\begin{equation}\label{zbound}
\mathbf{P}\left(\tilde{Z}_i\geq k^{1/2}\right)\leq\frac{\mathbf{E}(Z^{1+\alpha'})}{k^{\alpha'/2}}.
\end{equation}
Since the expectation in \eqref{zbound} is finite for any $\alpha'\in (0,\alpha-1)$ (see \cite[Section 35]{GneKol}, for example), we fix an $\alpha'$ in this range to obtain a polynomial bound for the relevant probability. For the second term of (\ref{twoterms}), we first condition on $\tilde{Z}_i$ to obtain
\begin{eqnarray*}
\mathbf{P}\left(\max_{j=1,\dots,\tilde{Z}_i-1}\#\mathcal{T}_{ij}\geq k^{1/2}-1\right)
&=&1-\mathbf{E}\left(\left(1-\mathbf{P}\left(\#\mathcal{T}\geq k^{1/2}-1\right)\right)^{\tilde{Z}-1}\right)\\
&=&1-\mathbf{E}\left(Z\left(1-\mathbf{P}\left(\#\mathcal{T}\geq k^{1/2}-1\right)\right)^{{Z}-1}\right)\\
&=&1-f'\left(1-\mathbf{P}\left(\#\mathcal{T}\geq k^{1/2}-1\right)\right).
\end{eqnarray*}
From the proof of Lemma \ref{bigtrap}, we know that $f'(1-x)\sim 1-\alpha x^{\alpha-1}L(x)$ as $x\rightarrow 0^+$, and so
\begin{equation}\label{tbound}
\mathbf{P}\left(\max_{j=1,\dots,\tilde{Z}_i-1}\#\mathcal{T}_{ij}\geq k^{1/2}-1\right)\sim \alpha \mathbf{P}\left(\#\mathcal{T}\geq k^{1/2}-1\right)^{\alpha-1}L\left(\mathbf{P}\left(\#\mathcal{T}\geq k^{1/2}-1\right)\right),
\end{equation}
as $k\rightarrow\infty$. To establish a bound for $\mathbf{P}(\#\mathcal{T} \geq k)$ that decays polynomially quickly, first note that $\mathbf{P}(\#\mathcal{T} = k)=k^{-1}\mathbf{P}(S_k=-1)$, where $(S_k)_{k\geq0}$ is a random walk on $\mathbb{Z}$ with step distribution $Z-1$ (see \cite{Dwass2}). Moreover, by the local limit theorem of \cite[Section 50]{GneKol}, it is the case that $\mathbf{P}(S_k=-1)\sim ca_k^{-1}$, where $a_k$ are the constants appearing in (\ref{stabledomain}). Since $a_k\sim k^{1/\alpha}\ell(k)$ for some slowly varying function $\ell$ (see \cite[Section 35]{GneKol}, for example), it follows that if $\alpha'\in (0,1/\alpha)$, then there exists a constant $c$ such that $\mathbf{P}(|\mathcal{T}| \geq k)\leq ck^{-\alpha'}$. Combining this estimate with (\ref{twoterms}), (\ref{zbound}) and (\ref{tbound}), we obtain that there exist constants $c$ and $\delta>0$ such that
\begin{equation}\label{vertexbound}
\mathbf{P}\left(\tilde{Z}_i\left[1+\max_{j=1,\dots,\tilde{Z}_i-1}\#\mathcal{T}_{ij}\right]\geq k\right)\leq ck^{-\delta}.
\end{equation}
Consequently, recalling \eqref{ebound},
\begin{eqnarray*}
\lefteqn{\mathbf{P}\left(\sum_{m\leq \Delta_{nT}}\mathbf{1}_{\{X_m\not\in\mathcal{B}\}}\geq\beta^{h_n(1+\varepsilon)}\right)}\\
&\leq &\mathbf{P}\left(\sum_{m\leq \Delta_{nT}}\mathbf{1}_{\{X_m\not\in\mathcal{B}\}}\geq n E^{\mathcal{T}^*}_{\rho}\left(\sum_{m\leq \Delta_{nT}}\mathbf{1}_{\{X_m\not\in\mathcal{B}\}}\right)\right)\\
&&\hspace{80pt}+\mathbf{P}\left(E^{\mathcal{T}^*}_{\rho}\left(\sum_{m\leq \Delta_{nT}}\mathbf{1}_{\{X_m\not\in\mathcal{B}\}}\right)\geq n^{-1}\beta^{h_n(1+\varepsilon)}\right)\\
&\leq & n^{-1}+ \mathbf{P}\left(\max_{i=0,\ldots,{nT-1}} \tilde{Z}_i\left[1+\max_{j=1,\dots,\tilde{Z}_i-1}\#\mathcal{T}_{ij}\right]\geq \frac{1}{c_\beta n^2T}\beta^{\varepsilon h_n}\right)\\
&\leq & n^{-1}+ nT\mathbf{P}\left(\tilde{Z}_i\left[1+\max_{j=1,\dots,\tilde{Z}_i-1}\#\mathcal{T}_{ij}\right]\geq \frac{1}{c_\beta n^2 T}\beta^{\varepsilon h_n}\right)\\
&\leq & n^{-1}+cn^{1+2\delta}\beta^{-\varepsilon\delta h_n},
\end{eqnarray*}
and this converges to 0 as $n\rightarrow\infty$.
\end{proof}
The result above means that, in establishing the distributional convergence of $\Delta_n$, we only have to consider the time the random walk $X$ spends in big leaves. In fact, as we will now show, the time spent close to the backbone in big leaves is also negligible. To this end, let us start by introducing some notation and formalising some terminology. First, we will write $y_{ij}$ for the deepest vertex in $\mathcal{T}_{ij}$; that is, the vertex that maximises the distance from the root $\rho_{ij}$. So that this notion is well-defined, if there is more than one vertex at the deepest level of $\mathcal{T}_{ij}$, we choose $y_{ij}$ to be the first in the usual lexicographical ordering of $\mathcal{T}_{ij}$, assuming that the offspring of each vertex have been labelled according to birth order. If the tree $\mathcal{T}_{ij}$ has height greater than or equal to $h_n$, then for a fixed $\delta\in(0,1)$ it is possible to define a unique vertex on the path from $\rho_{ij}$ to $y_{ij}$ at level $h_n^\delta$ in $\mathcal{T}_{ij}$. We shall denote this vertex $x_{ij}$ and call it the `entrance' to the leaf $\mathcal{T}_{ij}$. When we say that the leaf $\mathcal{T}_{ij}$ has been visited deeply, we will mean that $X$ has hit $x_{ij}$. Moreover, by the `time spent in the lower part of a big leaves emanating from $\rho_i$', we will mean
\begin{equation}\label{rrr}
t_i'=\sum_{j=1}^{\tilde{Z}_i-1}\mathbf{1}_{\{h(\mathcal{T}_{ij})\geq h_n\}}\sum_{m=0}^{\infty}\mathbf{1}_{\{X_m\in\mathcal{T}_{ij}\backslash \mathcal{T}_{ij}(x_{ij})\}},\end{equation}
where $\mathcal{T}_{ij}(x_{ij})$ is the part of the tree $\mathcal{T}_{ij}$ descending from the entrance $x_{ij}$.

To control the random variables $(t_i')_{i\geq 0}$  (which are identically-distributed apart from $i=0$), we need to consider the structure of the trees $\mathcal{T}_{ij}':=\mathcal{T}_{ij}\backslash \mathcal{T}_{ij}(x_{ij})$, and for this, the construction of a Galton-Watson tree conditioned on its height given in \cite{GK} is helpful. In particular, in Section 2 of that article, the following algorithm is described. First, let $(\xi_n,\zeta_n)$, $n\geq 0$, be a sequence of independent pairs of random variables, with distribution given by
\[\mathbf{P}\left(\xi_{n+1}=j,\zeta_{n+1}=k\right)=c_np_k\left(1-q_n\right)^{j-1}\left(1-q_{n+1}\right)^{k-j},\]
(recall that $q_n=\mathbf{P}(Z_n>0)$ is the probability that the unconditioned branching process survives for at least $n$ generations) for $1\leq j\leq k$, where
\begin{equation}\label{c}
c_n:=\frac{\mathbf{P}\left(h(\mathcal{T})=n\right)}{\mathbf{P}\left(h(\mathcal{T})=n+1\right)}.
\end{equation}
Then, let $\tilde{\mathcal{T}}_0$ be a Galton-Watson tree of height $0$, i.e. consisting solely of a root vertex, and, to construct $\tilde{\mathcal{T}}_{n+1}$, $n\geq 0$:
\begin{itemize}
  \item let the first generation size of $\tilde{\mathcal{T}}_{n+1}$ be $\zeta_{n+1}$,
  \item let $\tilde{\mathcal{T}}_n$ be the subtree founded by the $\xi_{n+1}$th first generation particle of $\tilde{\mathcal{T}}_{n+1}$,
  \item attach independent Galton-Watson trees conditioned on having height strictly less than $n$ to the $\xi_{n+1}-1$ siblings to the left of the distinguished first generation particle,
  \item attach independent Galton-Watson trees conditioned on height strictly less that $n+1$ to the $\zeta_{n+1}-\xi_{n+1}$ siblings to the right of the distinguished first generation particle.
\end{itemize}
It is shown in \cite{GK} that the tree $\tilde{\mathcal{T}}_n$ that results from this procedure has the same probabilistic structure as $\mathcal{T}$ conditioned to have height exactly equal to $n$. Before considering the implications of this result for the times $(t_i')_{i\geq 0}$, we derive the asymptotics of the constants $(c_n)_{n\geq 1}$ in our setting.

{\lem\label{clem} Let $\alpha\in(1,2]$. The constants $(c_n)_{n\geq 1}$, as defined at (\ref{c}), satisfy
\[c_n\sim 1+\frac{\alpha}{(\alpha-1)n},\]
as $n\rightarrow\infty$.}
\begin{proof} First note that $\mathbf{P}(h(\mathcal{T})=n)=q_n-q_{n+1}$. Moreover, if $f^{(n)}$ is the $n$-fold iteration of the generating function $f$, then we can write $q_n=1-f^{(n)}(0)$. It follows that
\[c_n=\frac{f^{(n+1)}(0)-f^{(n)}(0)}{f^{(n+2)}(0)-f^{ (n+1)}(0)}=\frac{f(1-q_n)-1+q_n}{f(1-q_{n+1})-1+q_{n+1}}=\frac{q_n^\alpha L(q_n)}{q_{n+1}^\alpha L(q_{n+1})}=\frac{q_n}{q_{n+1}}\times \frac{q_n^{\alpha-1} L(q_n)}{q_{n+1}^{\alpha-1} L(q_{n+1})},\]
where we have applied (\ref{genfun}) to deduce the third equality. Now, by (\ref{probdecay}), the second term on the right-hand side satisfies
\begin{equation}\label{qt1}
\frac{q_n^{\alpha-1} L(q_n)}{q_{n+1}^{\alpha-1} L(q_{n+1})}\sim \frac{n+1}{n}=1+\frac{1}{n}.
\end{equation}
For the first term, again applying (\ref{genfun}) and (\ref{probdecay}), it is the case that
\begin{equation}\label{qt2}
\frac{q_{n}}{q_{n+1}}=\frac{q_n}{1-f(1-q_n)}=\frac{1}{1-q_n^{\alpha-1}L(q_n)}\sim 1+\frac{1}{(\alpha-1)n}.
\end{equation}
Multiplying the right-hand sides of (\ref{qt1}) and (\ref{qt2}) yields the result. \end{proof}

{\lem \label{311} Let $\alpha\in(1,2]$ and $T\in(0,\infty)$. As $n\rightarrow\infty$,
\[\mathbb{P}_\rho\left(\sum_{i=0}^{nT-1}t_i'\geq \beta^{h_n}\right)\rightarrow 0.\]}
\begin{proof} Our first aim will be to show that
\begin{equation}\label{aim}
\mathbf{P}\left(h(\mathcal{T}_{ij}')>x\:\vline\:h(\mathcal{T}_{ij})\geq h_n\right)
\leq h_n^\delta \left(1-\left(\inf_{h\geq h_n-h_n^\delta}c_{h}\right)f'\left(1-q_{x-h_n^\delta}\right)\right),
\end{equation}
for $h_n^\delta<x<h_n$. Fix an $x$ in this range, and suppose for the moment that $h(\mathcal{T}_{ij})=h\geq h_n$, so that $x_{ij}$ is defined. Denote the path from $\rho_{ij}$ to $x_{ij}$ by $\rho_{ij}=w_0,w_1,\dots,w_{h_n^\delta}=x_{ij}$. Now, remove the edges $\{w_{l-1},w_{l}\}$, $l=1,\dots, h_n^\delta$ from $\mathcal{T}_{ij}'$, and denote by $\mathcal{T}_{ijl}$ the connected component containing $w_l$, so that $\mathcal{T}_{ij}'$ (minus the relevant edges) is the disjoint union of $\mathcal{T}_{ijl}$ over $l=0,\dots,h_n^\delta-1$. From the procedure for constructing a Galton-Watson tree conditioned on its height described before Lemma \ref{clem}, we deduce
\begin{eqnarray*}
\mathbf{P}\left(h(\mathcal{T}_{ij}')>x\:\vline\:h(\mathcal{T}_{ij})=h\right)&=&
\mathbf{P}\left(\max_{l=0,\dots,h_n^\delta-1}\left(h(\mathcal{T}_{ijl})+l\right)>x\:\vline\:h(\mathcal{T}_{ij})=h\right)\\
&\leq&\mathbf{P}\left(\max_{l=0,\dots,h_n^\delta-1} h(\mathcal{T}_{ijl})>x+1-h_n^\delta\:\vline\:h(\mathcal{T}_{ij})=h\right)\\
&=&1-\prod_{l=0}^{h_n^\delta-1}\left[1-\mathbf{P}\left(h(\mathcal{T}_{ijl})>x+1-h_n^\delta\:\vline\:h(\mathcal{T}_{ij})=h\right)\right].
\end{eqnarray*}
Moreover, if we suppose that $\mathcal{T}_{ij}$ conditioned on its height being equal to $h$ has been built from the random variables $(\xi_n,\zeta_n)$, $n\geq 0$, then we can write
\begin{eqnarray*}
\lefteqn{
\mathbf{P}\left(h(\mathcal{T}_{ijl})<x+1-h_n^\delta\:\vline\:h(\mathcal{T}_{ij})=h\right)}\\
&=&
\mathbf{E}\left(\mathbf{P}\left(h(\mathcal{T}_{ijl})<x+1-h_n^\delta\:\vline\:h(\mathcal{T}_{ij})=h, \xi_{h-l},\zeta_{h-l}\right)\right)\\
&=&\mathbf{E}\left(\mathbf{P}\left(h(\mathcal{T})<x-h_n^\delta\:\vline\:h(\mathcal{T})<h-l-1\right)^{\xi_{h-l}-1}
\vphantom{\mathbf{P}\left(h(\mathcal{T})<x-h_n^\delta\:\vline\:h(\mathcal{T})<h-l\right)^{\zeta_{h-l}-\xi_{h-l}}}\right.\\
&&\hspace{50pt}\times\left.\vphantom{\mathbf{P}\left(h(\mathcal{T})<x-h_n^\delta\:\vline\:h(\mathcal{T})<h-l-1\right)^{\xi_{h-l}-1}}
\mathbf{P}\left(h(\mathcal{T})<x-h_n^\delta\:\vline\:h(\mathcal{T})<h-l\right)^{\zeta_{h-l}-\xi_{h-l}}\right)\\
&=&\sum_{1\leq j\leq k}c_{h-l-1}p_k\mathbf{P}\left(h(\mathcal{T})<h-l-1\right)^{j-1}\mathbf{P}\left(h(\mathcal{T})<x-h_n^\delta\:\vline\:h(\mathcal{T})<h-l-1\right)^{j-1}\\
&&\hspace{80pt}\times\mathbf{P}\left(h(\mathcal{T})<h-l\right)^{k-j}\mathbf{P}\left(h(\mathcal{T})<x-h_n^\delta\:\vline\:h(\mathcal{T})<h-l\right)^{k-j}\\
&=&\sum_{1\leq j\leq k}c_{h-l-1}p_k\mathbf{P}\left(h(\mathcal{T})<x-h_n^\delta\right)^{j-1}\mathbf{P}\left(h(\mathcal{T})<x-h_n^\delta\right)^{k-j}\\
&=&\sum_{k\geq 1}c_{h-l-1} k p_k\mathbf{P}\left(h(\mathcal{T})<x-h_n^\delta\right)^{k-1}\\
&=&c_{h-l-1}f'\left(1-q_{x-h_n^\delta}\right).
\end{eqnarray*}
Thus, combining these deductions, we obtain
\begin{eqnarray*}
\mathbf{P}\left(h(\mathcal{T}_{ij}')>x\:\vline\:h(\mathcal{T}_{ij})=h\right)
&\leq & 1-\left(c_{h-l-1}f'\left(1-q_{x-h_n^\delta}\right)\right)^{h_n^\delta}\\
&\leq & h_n^\delta \left(1-\inf_{h'\geq h_n-h_n^\delta}c_{h'}f'\left(1-q_{x-h_n^\delta}\right)\right),
\end{eqnarray*}
and, since this bound is independent of $h\geq h_n$, the bound at (\ref{aim}) follows.

Now, by arguing similarly to (\ref{ebound}), it is possible to check that
\[E_{\rho}^{\mathcal{T}^*}\left(\sum_{i=0}^{nT-1}t_i'\right)\leq c_\beta\sum_{i=0}^{nT-1}\sum_{j=1}^{\tilde{Z}_i-1}\mathbf{1}_{\{h(\mathcal{T}_{ij})\geq h_n\}}\beta^{h(\mathcal{T}_{ij}')}\#\mathcal{T}_{ij},\]
where $c_\beta$ is a constant depending only upon $\beta$. Thus, following the end of the proof of Lemma \ref{smallleaf},
\begin{eqnarray*}
\lefteqn{\mathbb{P}_\rho\left(\sum_{i=0}^{nT-1}t_i'\geq \beta^{h_n}\right)}\\
&\leq &n^{-1}+nT\mathbf{P}\left(T c_\beta\sum_{j=1}^{\tilde{Z}_i-1}\mathbf{1}_{\{h(\mathcal{T}_{ij})\geq h_n\}}\beta^{h(\mathcal{T}_{ij}')}\#\mathcal{T}_{ij}\geq n^{-2}\beta^{h_n}\right)\\
&\leq & n^{-1}+nT\mathbf{P}\left(T c_\beta\tilde{Z}_i\max_{j=1,\dots,\tilde{Z}_i-1}\#\mathcal{T}_{ij}\geq n^{-2}\beta^{h_n/2}\right)+nT\mathbf{P}\left(\max_{\substack{j=1,\dots,\tilde{Z}_i-1:\\h(\mathcal{T}_{ij})\geq h_n}}{h(\mathcal{T}_{ij}')}\geq {h_n/2}\right).
\end{eqnarray*}
Clearly the first term decays to zero, and, by applying (\ref{vertexbound}), so does the second term. To deal with the third term, observe that, under the convention that $h(\mathcal{T}_{ij}')=0$ for $j$ such that $h(\mathcal{T}_{ij})<h_n$,
\begin{eqnarray*}
\mathbf{P}\left(\max_{\substack{j=1,\dots,\tilde{Z}_i-1:\\h(\mathcal{T}_{ij})\geq h_n}}{h(\mathcal{T}_{ij}')}\geq {h_n/2}\right)& = &
1-\mathbf{E}\left(\left(1-\mathbf{P}\left(h(\mathcal{T}_{ij}')\geq h_n/2\right)\right)^{\tilde{Z}-1}\right)\\
&=&1-f'\left(1-\mathbf{P}\left(h(\mathcal{T}_{ij}')\geq h_n/2\right)\right)\\
&\sim&\alpha\mathbf{P}\left(h(\mathcal{T}_{ij}')\geq h_n/2\right)^{\alpha-1}L\left(\mathbf{P}\left(h(\mathcal{T}_{ij}')\geq h_n/2\right)\right)\\
&\sim& \alpha q_{h_n}^{\alpha-1} \mathbf{P}\left(h(\mathcal{T}_{ij}')>h_n/2\:\vline\:h(\mathcal{T}_{ij})\geq h_n\right)^{\alpha-1}\\
&&\hspace{50pt}L\left(q_{h_n} \mathbf{P}\left(h(\mathcal{T}_{ij}')>h_n/2\:\vline\:h(\mathcal{T}_{ij})\geq h_n\right)\right),
\end{eqnarray*}
where we have used that $f'(1-x)\sim 1-\alpha x^{\alpha-1}L(x)$ as $x\rightarrow 0^+$, which we first recalled in the proof of Lemma \ref{bigtrap}, and (\ref{probdecay}) again. Since the representation theorem for slowly varying functions (\cite[Theorem 1.2]{sen}, for example) implies that, for any $\varepsilon>0$,
\[L\left(q_{h_n} \mathbf{P}\left(h(\mathcal{T}_{ij}')>h_n/2\:\vline\:h(\mathcal{T}_{ij})\geq h_n\right)\right)\leq \mathbf{P}\left(h(\mathcal{T}_{ij}')>h_n/2\:\vline\:h(\mathcal{T}_{ij})\geq h_n\right)^{-\varepsilon} L(q_{h_n}),\]
for large $n$, it follows that $\mathbf{P}(\max_{{j=1,\dots,\tilde{Z}_i-1:\:h(\mathcal{T}_{ij})\geq h_n}}\max_{j\in B_i}{h(\mathcal{T}_{ij}')}\geq {h_n/2})$ is asymptotically less than
\[\alpha\mathbf{P}\left(h(\mathcal{T}_{ij}')>h_n/2\:\vline\:h(\mathcal{T}_{ij})\geq h_n\right)^{\alpha-1-\varepsilon}q_{h_n}^{\alpha-1}L(q_{h_n})\sim \frac{\alpha\mathbf{P}\left(h(\mathcal{T}_{ij}')>h_n/2\:\vline\:h(\mathcal{T}_{ij})\geq h_n\right)^{\alpha-1-\varepsilon}}{(\alpha-1)h_n}.\]
Finally, setting $x=h_n/2$ in (\ref{aim}) and applying Lemma \ref{clem} yields
\begin{eqnarray*}
\mathbf{P}\left(h(\mathcal{T}_{ij}')>h_n/2\:\vline\:h(\mathcal{T}_{ij})\geq h_n\right)
&\leq& h_n^\delta \left(1-\inf_{h\geq h_n-h_n^\delta}c_hf'\left(1-q_{ 2^{-1}h_n-h_n^\delta}\right)\right)\\
&\sim & \alpha h_n^\delta q_{2^{-1} h_n-h_n^\delta}^{\alpha-1}L\left(q_{2^{-1} h_n-h_n^\delta}\right)\\
&\sim& c h_n^{\delta-1},
\end{eqnarray*}
for a suitable choice of constant $c$, and so, by adjusting $c$ as necessary, we obtain that, for large $n$,
\[nT\mathbf{P}\left(\max_{\substack{j=1,\dots,\tilde{Z}_i-1:\\h(\mathcal{T}_{ij})\geq h_n}}{h(\mathcal{T}_{ij}')}\geq {h_n/2}\right)\leq cnh_n^{-1}  h_n^{(\delta-1)(\alpha-1-\varepsilon)}.\]
Since this upper bound converges to 0 for any $\varepsilon<\alpha-1$, this completes the proof.
\end{proof}

In deriving tail asymptotics for the time $X$ spends in the big leaves emanating from a particular backbone vertex, it will be useful to have information about the set of big leaves that the biased random walk visits deeply before it escapes along the backbone, and the next two lemmas provide this. For their statement, we define the index set of big leaves emanating from $\rho_i$ by
\[B_i:=\left\{j=1,\dots,\tilde{Z}_i-1:\:h(\mathcal{T}_{ij})\geq h_n\right\}\]
and the subset of those that are visited deeply by $X$ before it escapes a certain distance along the backbone by
\[V_i:=\left\{j\in B_i:\:\tau_{x_{ij}}<\tau_{z_i}\right\},\]
where $z_i:=\rho_{i+1+h_n^\delta}$.

{\lem\label{vlem} Let $\alpha\in(1,2]$ and $i\geq0$. For any $A\subseteq B_i$, we have
\[P^{\mathcal{T}^*}_{\rho_i}\left(V_i=A\right)=\frac{1}{1+\#B_i}\binom{\#B_i}{\#A}^{-1},\]
where $\#B_i$, $\#A$ represents the cardinality of $B_i$, $A$, respectively.}
\begin{proof} The lemma readily follows from the symmetry of the situation, which implies that, starting from $\rho_i$, the biased random walk $X$ is equally likely to visit any one of $x_{ij}$, $j\in B_i$ and $z_i$ first.
\end{proof}

Although the above lemma might seem simple, it allows us to deduce the distributional tail behaviour of the greatest height of a big leaf at a particular backbone vertex visited by the biased random walk $X$. Note that we continue to use the notation $q_n=\mathbf{P}(Z_n>0)$.

{\lem \label{hseen} Let $\alpha\in(1,2]$ and $i\geq 0$. For $x\geq h_n$,
\[\mathbb{P}_\rho\left(\max_{j\in V_i}h(\mathcal{T}_{ij})\geq x\right)=q_x^{\alpha-1}L\left(q_{x}\right).\]}
\begin{proof} Let $x\geq h_n$. By definition, we have that
\[\mathbb{P}_\rho\left(\max_{j\in V_i}h(\mathcal{T}_{ij})< x\right)=
\mathbf{E}\left({E}_\rho^{\mathcal{T}^*}\left(\prod_{j\in V_i}\mathbf{1}_{\{h(\mathcal{T}_{ij})<x\}}\right)\right),\]
and decomposing the inner expectation over the possible values of $V_i$ yields
\[{E}_\rho^{\mathcal{T}^*}\left(\prod_{j\in V_i}\mathbf{1}_{\{h(\mathcal{T}_{ij})<x\}}\right)
=\sum_{A\subseteq B_i}{E}_\rho^{\mathcal{T}^*}\left(\mathbf{1}_{\{V_i=A\}}\prod_{j\in A}\mathbf{1}_{\{h(\mathcal{T}_{ij})<x\}}\right).\]
Since $\prod_{j\in A}\mathbf{1}_{\{h(\mathcal{T}_{ij})<x\}}$ is a measurable function of $\mathcal{T}^*$, this can be rewritten as
\begin{eqnarray*}
{E}_\rho^{\mathcal{T}^*}\left(\prod_{j\in V_i}\mathbf{1}_{\{h(\mathcal{T}_{ij})<x\}}\right)
&=&\sum_{A\subseteq B_i} {P}_\rho^{\mathcal{T}^*}\left(V_i=A\right)\prod_{j\in A}\mathbf{1}_{\{h(\mathcal{T}_{ij})<x\}}\\
&=&\sum_{A\subseteq B_i} \frac{1}{\#B_i+1}\binom{\#B_i}{\#A}^{-1}\prod_{j\in A}\mathbf{1}_{\{h(\mathcal{T}_{ij})<x\}},
\end{eqnarray*}
where the second equality is an application of Lemma \ref{vlem}. Now, since
\[\mathbf{P}\left(\prod_{j\in A}\mathbf{1}_{\{h(\mathcal{T}_{ij})<x\}}\:\vline\:B_i\right)=\mathbf{P}\left(h(\mathcal{T})< x\:\vline\:h(\mathcal{T})\geq h_n\right)^{\#A}= \left(1-\frac{q_x}{q_{h_n}}\right)^{\#A},\]
for every $A\subseteq B_i$, it follows that
\begin{eqnarray*}
\mathbb{P}_\rho\left(\max_{j\in V_i}h(\mathcal{T}_{ij})< x\right)&=&
\mathbf{E}\left(\sum_{A\subseteq B_i} \frac{1}{\#B_i+1}\binom{\#B_i}{\#A}^{-1}
\mathbf{P}\left(\prod_{j\in A}\mathbf{1}_{\{h(\mathcal{T}_{ij})<x\}}\:\vline\:B_i\right)\right)\\
&=&
\mathbf{E}\left(\sum_{A\subseteq B_i} \frac{1}{\#B_i+1}\binom{\#B_i}{\#A}^{-1}
\left(1-\frac{q_x}{q_{h_n}}\right)^{\#A}\right)\\
&=&\mathbf{E}\left(\frac{1}{\#B_i+1}\sum_{l=0}^{\#B_i}\left(1-\frac{q_x}{q_{h_n}}\right)^{l}\right)\\
&=&\mathbf{E}\left(\frac{q_{h_n}}{(\#B_i+1)q_x}\left(1-\left(1-\frac{q_x}{q_{h_n}}\right)^{\#B_i+1}\right)\right).
\end{eqnarray*}
To continue, observe that, conditional on $\tilde{Z}_i$, $\#B_i$ is binomially distributed with parameters $\tilde{Z}_i-1$ and $q_{h_n}$. Consequently, the probability we are trying to compute is equal to
\begin{equation}\label{t0}
\mathbf{E}\left(\sum_{l=0}^{\tilde{Z}-1}\binom{\tilde{Z}-1}{l}q_{h_n}^{l}\left(1-q_{h_n}\right)^{\tilde{Z}-1-l}\frac{q_{h_n}}{(l+1)q_x}\left(1-\left(1-\frac{q_x}{q_{h_n}}\right)^{l+1}\right)\right).
\end{equation}
We break this into two terms. Firstly,
\begin{eqnarray}
\lefteqn{\mathbf{E}\left(\sum_{l=0}^{\tilde{Z}-1}\binom{\tilde{Z}-1}{l}q_{h_n}^{l}\left(1-q_{h_n}\right)^{\tilde{Z}-1-l}\frac{q_{h_n}}{(l+1)q_x}\right)}\nonumber\\
&=&q_x^{-1}\mathbf{E}\left(\sum_{l=0}^{{Z}-1}\binom{Z}{l+1}q_{h_n}^{l+1}\left(1-q_{h_n}\right)^{{Z}-1-l}\mathbf{1}_{\{Z\geq 1\}}\right)\nonumber\\
&=&q_x^{-1}\mathbf{E}\left(1-\left(1-q_{h_n}\right)^Z\right)\nonumber\\
&=&q_x^{-1}\left(1-f\left(1-q_{h_n}\right)\right).\label{t1}
\end{eqnarray}
Secondly,
\begin{eqnarray}
\lefteqn{\mathbf{E}\left(\sum_{l=0}^{\tilde{Z}-1}\binom{\tilde{Z}-1}{l}q_{h_n}^{l}\left(1-q_{h_n}\right)^{\tilde{Z}-1-l}\frac{q_{h_n}}{(l+1)q_x}\left(1-\frac{q_x}{q_{h_n}}\right)^{l+1}\right)}\nonumber\\
&=&q_x^{-1}\mathbf{E}\left(\sum_{l=0}^{{Z}-1}\binom{Z}{l+1}\left(q_{h_n}-q_x\right)^{l+1}\left(1-q_{h_n}\right)^{{Z}-1-l}\mathbf{1}_{\{Z\geq 1\}}\right)\nonumber\\
&=&q_x^{-1}\left(f\left(1-q_x\right)-f\left(1-q_{h_n}\right)\right).\label{t2}
\end{eqnarray}
Since taking the difference between (\ref{t1}) and (\ref{t2}) gives us (\ref{t0}), we have thus proved that
\begin{eqnarray*}
\mathbb{P}_\rho\left(\max_{j\in V_i}h(\mathcal{T}_{ij})< x\right)&=&
q_x^{-1}\left(1-f\left(1-q_x\right)\right)\\
&=&q_x^{-1}\left(q_x-q_x^\alpha L\left(q_x\right)\right)\\
&=&1-q_x^{\alpha-1}L\left(q_{x}\right),
\end{eqnarray*}
where the second equality is a consequence of (\ref{genfun}), and the lemma follows.
\end{proof}

With these preparations in place, we are now ready to study the asymptotic tail behaviour of
\[t_i:=\sum_{j=1}^{\tilde{Z}_i-1}\mathbf{1}_{\{h(\mathcal{T}_{ij})\geq h_n\}}\sum_{m=0}^{\tau_{z_i}}\mathbf{1}_{\{X_m\in \mathcal{T}_{ij}(x_{ij})\}},\]
which can be interpreted as the length of time the $X$ spends deep inside leaves emanating from $\rho_i$ before escaping along the backbone. The next lemma gives an upper tail bound for this random variable.

{\lem \label{deeptime} Let $\alpha\in(1,2]$ and $\varepsilon>0$. There exists a constant $c_{\beta,\varepsilon}$ such that, for any $i\geq 0$ and $x$ satisfying $\ln x\geq c_{\beta,\varepsilon}h_n$,
\[\mathbb{P}_\rho\left(t_i\geq x\right)\leq  \frac{(1+\varepsilon)\ln \beta}{(\alpha-1)\ln x}.\]}
\begin{proof} First note that, by applying the commute time identity for random walks (for example, \cite[Proposition 10.6]{LPW}), we have that
\[E^{\{\rho_i\}\cup\mathcal{T}_{ij}}_{x_{ij}}\tau_{\rho_i}+E^{\{\rho_i\}\cup\mathcal{T}_{ij}}_{\rho_i}\tau_{x_{ij}}=\left(\beta^{-i}+\dots+\beta^{-(i+h_n^\delta)}\right)\left(2\beta^{i}+\sum_{\substack{x,y\in \mathcal{T}_{ij}:\\ x\sim y}}c(x,y)\right),\]
where $E^{\{\rho_i\}\cup\mathcal{T}_{ij}}_\cdot$ refers to the random walk on the tree $\mathcal{T}_{ij}$ extended by adding the vertex $\rho_i$ and the edge $\{\rho_i,\rho_{ij}\}$. Since $E^{\{\rho_i\}\cup\mathcal{T}_{ij}}_{x_{ij}}\tau_{\rho_i}=E^{\mathcal{T}^*}_{x_{ij}}\tau_{\rho_i}$, it follows that
\[E^{\mathcal{T}^*}_{x_{ij}}\tau_{\rho_i}\leq \frac{\beta}{\beta-1}\left(2+2\beta^{h(\mathcal{T}_{ij})}\#\mathcal{T}_{ij}\right).\]
Thus, since the random walk $X$ spends no time in $\mathcal{T}_{ij}(x_{ij})$ if $j\not\in V_i$, we can bound the quenched expectation of $t_i$ conditional on $V_i$ as follows:
\begin{eqnarray}
E^{\mathcal{T}^*}_\rho \left( t_i\:\vline\:V_i\right)&=&\sum_{j=1}^{\tilde{Z}_i-1}\sum_{m=0}^{\infty}P^{\mathcal{T}^*}_\rho\left(X_m\in \mathcal{T}_{ij}(x_{ij}),\:m\leq \tau_{z_i}\:\vline\:V_i\right)\nonumber\\
&\leq& \sum_{j\in V_i} E^{\mathcal{T}^*}_{x_{ij}}\left(\tau_{\rho_i}\right)
E^{\mathcal{T}^*}_{\rho}\left(\upsilon_{ij}\:|\:V_i\right),  \label{b2} \\
&\leq & c_{\beta}\sum_{j\in V_i}\left(1+\#\mathcal{T}_{ij}\right)\beta^{h(\mathcal{T}_{ij})}E^{\mathcal{T}^*}_{\rho}\left(\upsilon_{ij}\:|\:V_i\right)\label{b3},
\end{eqnarray}
where $\upsilon_{ij}$ is the number of passages $X$ makes from $\rho_i$ to $x_{ij}$ before it hits $z_i$, and the inequality at (\ref{b2}) is obtained by an application of the strong Markov property (that holds with respect to the unconditioned law). Now, $\upsilon_{ij}$ is clearly bounded above by the total number of visits to $\rho_i$, $N(\rho_i)$ say, and, by symmetry, this latter random variable satisfies $E^{\mathcal{T}^*}_{\rho}(N(\rho_i)\:|\:V_i) =E^{\mathcal{T}^*}_{\rho}(N(\rho_i)\:|\:\#V_i)$. Consequently, we deduce that
\begin{eqnarray*}
E^{\mathcal{T}^*}_{\rho}\left(\upsilon_{ij}\:|\:V_i\right)&\leq&
\sum_{k=0}^{\#B_i}\mathbf{1}_{\{\#V_i=k\}}E^{\mathcal{T}^*}_{\rho}\left(N(\rho_i)\:|\:\#V_i=k\right)\\
&\leq&\sum_{k=0}^{\#B_i}\mathbf{1}_{\{\#V_i=k\}}\frac{E^{\mathcal{T}^*}_{\rho}\left(N(\rho_i)\right)}{P^{\mathcal{T}^*}_{\rho}\left(\#V_i=k\right)}\\
&\leq&c_\beta\sum_{k=0}^{\#B_i}\mathbf{1}_{\{\#V_i=k\}}(\#B_i+1)\tilde{Z}_i\\
&\leq&c_\beta \tilde{Z}_i^2,
\end{eqnarray*}
where we have applied Lemma \ref{vlem} and the argument at (\ref{est2}) to deduce $P^{\mathcal{T}^*}_{\rho}(\#V_i=k)^{-1}=\#B_i+1$ and $E^{\mathcal{T}^*}_{\rho}(N(\rho_i))\leq c_{\beta}\tilde{Z}_i$, respectively.
Applying the above bound in combination with (\ref{b3}) yields
\[E^{\mathcal{T}^*}_\rho \left( t_i\:\vline\:V_i\right)\leq  c_\beta \tilde{Z}_i^3 \left(1+\max_{j=1,\dots,\tilde{Z}_i-1}\#\mathcal{T}_{ij}\right)\beta^{\max_{j\in V_i}h(\mathcal{T}_{ij})}.\]
Thus, for $\eta\in(0,\frac{1}{2})$, we can conclude
\begin{eqnarray*}
\mathbb{P}_\rho\left(t_i\geq x\right)&\leq & \mathbb{P}_\rho\left(t_i\geq x^\eta E^{\mathcal{T}^*}_\rho \left( t_i\:\vline\:V_i\right)\right)+\mathbb{P}_\rho\left(E^{\mathcal{T}^*}_\rho \left( t_i\:\vline\:V_i\right)\geq x^{1-\eta} \right)\\
&\leq &x^{-\eta}+\mathbb{P}_\rho\left(c_\beta \tilde{Z}_i^3 \left(1+\max_{j=1,\dots,\tilde{Z}_i-1}\#\mathcal{T}_{ij}\right)\beta^{\max_{j\in V_i}h(\mathcal{T}_{ij})}\geq x^{1-\eta}\right)\\
&\leq &x^{-\eta}+\mathbf{P}\left(c_\beta \tilde{Z}_i^3 \left(1+\max_{j=1,\dots,\tilde{Z}_i-1}\#\mathcal{T}_{ij}\right)^2\geq x^\eta\right)
+\mathbb{P}_\rho\left(\beta^{\max_{j\in V_i}h(\mathcal{T}_{ij})}\geq x^{1-2\eta}\right)\\
&\leq & x^{-\eta}+c_\beta x^{-\eta\delta/3}+q_{(1-2\eta)\ln x/\ln \beta}^{\alpha-1}L\left(q_{(1-2\eta)\ln x/\ln \beta}\right)
\end{eqnarray*}
for $(1-2\eta)\ln x/\ln \beta\geq h_n$, where the value of $c_\beta$ has been updated from above and the constant $\delta$ is the one appearing in (\ref{vertexbound}). We have also applied Lemma \ref{hseen} in obtaining the final bound. Finally, (\ref{probdecay}) allows us to deduce from this that, as long as $(1-2\eta)\ln x/\ln \beta$ is sufficiently large, it holds that
\[\mathbb{P}_\rho\left(t_i\geq x\right)\leq  \frac{(1+\eta)\ln \beta}{(1-2\eta)(\alpha-1)\ln x}.\]
The result follows.
\end{proof}

We can also prove a lower bound for the distributional tail of $t_i$ that matches the upper bound proved above. Similarly to a proof strategy followed in \cite{AFGH}, a key step in doing this is obtaining a concentration result to show that the time spent in a leaf visited deeply by the process $X$ will be on the same scale as its expectation.

{\lem \label{deeptimelower} Let $\alpha\in(1,2]$ and $\varepsilon>0$. There exist constants $n_0$ and $c_{\beta,\varepsilon}$ such that, for any $i\geq 0$, $n\geq n_0$ and $x$ satisfying $c_{\beta,\varepsilon}h_n \leq \ln x\leq n^2$,
\[\mathbb{P}_\rho\left(t_i\geq x\right)\geq  \frac{(1-\varepsilon)\ln \beta}{(\alpha-1)\ln x}.\]}
\begin{proof} Our first goal is to derive an estimate on the lower tail of the time that $X$ spends in a big leaf $\mathcal{T}_{ij}$ before hitting $\rho_i$, given that it starts at the entrance vertex $x_{ij}$. To this end, we start by noting that under $P^{\mathcal{T}^*}_{x_{ij}}$ and conditional on the number of returns that the random walk $X$ makes to $\mathcal{T}_{ij}(x_{ij})$ before hitting $\rho_i$, i.e.
\[\upsilon'_{ij}:=\#\left\{m\leq \tau_{\rho_i}:\:X_{m-1}=x_{ij}',\:X_m=x_{ij}\right\},\]
where $x_{ij}'$ denotes the parent of $x_{ij}$, the random variable $\Sigma:=\sum_{m=0}^{\tau_{\rho_i}}\mathbf{1}_{\{X_m\in\mathcal{T}_{ij}(x_{ij})\}}$
is distributed as $\upsilon'_{ij}+1$ independent copies of a random variable whose law is equal to that of $\tau_{x_{ij}'}$ under $P^{\mathcal{T}^*}_{x_{ij}}$. (This is a simple application of the strong Markov property.) In particular, we have that
\[E_{x_{ij}}^{\mathcal{T}^*}\left(\Sigma\:\vline\:\upsilon'_{ij}\right)=
\left(1+\upsilon'_{ij}\right)E_{x_{ij}}^{\mathcal{T}^*}\tau_{x_{ij}'},\]
and also
\[{\rm Var}_{x_{ij}}^{\mathcal{T}^*}\left(\Sigma\:\vline\:\upsilon'_{ij}\right)=\left(1+\upsilon'_{ij}\right)
{\rm Var}_{x_{ij}}^{\mathcal{T}^*}\left(\tau_{x_{ij}'}\right).\]
To control the right-hand sides of these quantities, we will apply the following moment bounds:
\[E_{x_{ij}}^{\mathcal{T}^*}\tau_{x_{ij}'}=1+\beta^{-(i+h_n^\delta)}\sum_{\substack{x,y\in\mathcal{T}_{ij}(x_{ij}):\\x\sim y}}c(x,y)\geq \beta^{h(\mathcal{T}_{ij})-h_n^\delta}\]
and
\[E_{x_{ij}}^{\mathcal{T}^*}\left(\tau_{x_{ij}'}^2\right)\leq \frac{2\beta}{\beta-1}E_{x_{ij}}^{\mathcal{T}^*}\left(\tau_{x_{ij}'}\right)^2,\]
where the first moment lower bound is obtained by applying a formula similar to (\ref{treeexit}), and the second moment upper bound is an adaptation of a result derived in the proof of \cite[Lemma 9.1]{AFGH}. As for the distribution of $\upsilon'_{ij}$ under $P^{\mathcal{T}^*}_{x_{ij}}$, it is clear this is geometric, with parameter given by
\[P_{x_{ij}'}^{\mathcal{T}^*}\left(\tau_{\rho_i}<\tau_{x_{ij}}\right)=\frac{\left(1+\beta^{-1}+\dots+\beta^{-h_n^\delta+1}\right)^{-1}}{\beta^{h_n^\delta}+\left(1+\beta^{-1}+\dots+\beta^{-h_n^\delta+1}\right)^{-1}}=\frac{\beta^{h_n^\delta}-\beta^{h_n^\delta-1}}{\beta^{2h_n^\delta}-\beta^{h_n^\delta-1}},\]
from which it follows that
\[E_{x_{ij}'}^{\mathcal{T}^*}\left(\upsilon'_{ij}+1\right)=\frac{\beta^{h_n^\delta}-\beta^{-1}}{1-\beta^{-1}}\geq \beta^{h_n^\delta}-\beta^{-1}\geq \frac{\beta^{h_n^\delta}}{2},\]
for $n\geq n_0$, where $n_0$ is a deterministic constant. Putting the above observations together, we deduce that, for $n\geq n_0$ and $\varepsilon>0$,
\begin{eqnarray*}
\lefteqn{P^{\mathcal{T}^*}_{x_{ij}}\left(\Sigma\leq \frac{\varepsilon}{4}\beta^{h(\mathcal{T}_{ij})}\right)}\\
&\leq&P^{\mathcal{T}^*}_{x_{ij}}\left(\Sigma\leq \frac{\varepsilon}{4}\beta^{h(\mathcal{T}_{ij})},\:\upsilon'_{ij}+1\geq\varepsilon E^{\mathcal{T}^*}_{x_{ij}} (\upsilon'_{ij}+1)\right)+
P^{\mathcal{T}^*}_{x_{ij}}\left(\upsilon'_{ij}+1<\varepsilon E^{\mathcal{T}^*}_{x_{ij}} (\upsilon'_{ij}+1)\right)\\
&\leq &P^{\mathcal{T}^*}_{x_{ij}}\left(\Sigma\leq \frac{1}{2}\left(\upsilon'_{ij}+1\right)E_{x_{ij}}^{\mathcal{T}^*}\tau_{x_{ij}'}\right)+1-\left(1-\frac{\beta^{h_n^{\delta}}-\beta^{h_n^{\delta}-1}}{\beta^{2h_n^{\delta}}-\beta^{h_n^{\delta}-1}}\right)^{\varepsilon E^{\mathcal{T}^*}_{x_{ij}} (\upsilon'_{ij}+1)} \\
&\leq&P^{\mathcal{T}^*}_{x_{ij}}\left(\left|
\Sigma-E_{x_{ij}}^{\mathcal{T}^*}\left(\Sigma\:\vline\:\upsilon'_{ij}\right)\right|\geq \frac{1}{2}E_{x_{ij}}^{\mathcal{T}^*}\left(\Sigma\:\vline\:\upsilon'_{ij}\right)\right)+\varepsilon\\
&\leq& E^{\mathcal{T}^*}_{x_{ij}}\left(\frac{4{\rm Var}_{x_{ij}}^{\mathcal{T}^*}\left(\Sigma\:\vline\:\upsilon'_{ij}\right)}
{E_{x_{ij}}^{\mathcal{T}^*}\left(\Sigma\:\vline\:\upsilon'_{ij}\right)^2}\right)+\varepsilon\\
&\leq &\frac{8\beta}{\beta-1}E^{\mathcal{T}^*}_{x_{ij}}\left(\frac{1}{\upsilon'_{ij}+1}\right)+\varepsilon\\
&\leq & c h_n^\delta\beta^{-h_n^\delta}+\varepsilon,
\end{eqnarray*}
where $c$ is a constant depending only on $\beta$ and $n_0$ (and not $\varepsilon$).

Now, if we suppose $j_0\in V_i$ is such that $h(\mathcal{T}_{ij_0})=\max_{j\in V_i}h(\mathcal{T}_{ij})$, then
\[P^{\mathcal{T}^*}_\rho\left(t_i\leq \frac{\varepsilon}{4}\beta^{\max_{j\in V_i}h(\mathcal{T}_{ij})}\:\vline\:V_i\right)\leq
P^{\mathcal{T}^*}_\rho\left(\tau_{x_{ij_0}\rightarrow \rho_i}\leq \frac{\varepsilon}{4}\beta^{h(\mathcal{T}_{ij_0})}\:\vline\:V_i\right),\]
where $\tau_{x_{ij_0}\rightarrow \rho_i}$ is the amount of time $X$ spends in $\mathcal{T}_{ij_0}(x_{ij_0})$ before $\inf\{m\geq \tau_{x_{ij_0}}:\:X_m=\rho_i\}$. By applying a strong Markov argument for the unconditioned law (cf. (\ref{b2})), yields that the law of $\tau_{x_{ij_0}\rightarrow \rho_i}$ under $P^{\mathcal{T}^*}_{\rho}(\cdot|\:V_i)$ is the same as that of $\Sigma$ (as defined above with $j=j_0$) under $P^{\mathcal{T}^*}_{x_{ij_0}}$, and thus the result of the previous paragraph implies that, for $n\geq n_0$ and $\varepsilon>0$,
\[P^{\mathcal{T}^*}_\rho\left(t_i\leq \frac{\varepsilon}{4}\beta^{\max_{j\in V_i}h(\mathcal{T}_{ij})}\:\vline\:V_i\right)\leq c h_n^\delta\beta^{-h_n^\delta}+\varepsilon.\]
Taking expectations with respect to $P^{\mathcal{T}^*}_\rho$ and $\mathbf{P}$ establishes that the same is true when $P^{\mathcal{T}^*}_{\rho}(\cdot|\:V_i)$ is replaced by the annealed law $\mathbb{P}_\rho$. Consequently, for any $n\geq n_0$, $\varepsilon>0$ and $\ln x\geq h_n\ln \beta$,
\begin{eqnarray*}
\mathbb{P}_\rho\left(t_i\leq x\right)&\leq &\mathbb{P}_\rho\left(t_i\leq \frac{\varepsilon}{4}\beta^{\max_{j\in V_i}h(\mathcal{T}_{ij})}\right)+\mathbb{P}_\rho\left(\max_{j\in V_i}h(\mathcal{T}_{ij})\leq
\frac{\ln (4x/\varepsilon)}{\ln \beta}\right)\\
&\leq & c h_n^\delta\beta^{-h_n^\delta}+\varepsilon + 1-q_{\frac{\ln (4x/\varepsilon)}{\ln \beta}}^{\alpha-1}L\left(q_{\frac{\ln (4x/\varepsilon)}{\ln \beta}}\right),
\end{eqnarray*}
where we have applied Lemma \ref{hseen} to deduce the second inequality. Finally, fix $\eta>0$. If we set $\varepsilon=1/(\ln x)^2$, then the second term is bounded above by $\eta/\ln x$ for any $\ln x\geq \eta^{-1}$. With this choice of $\varepsilon$, by (\ref{probdecay}), the fourth term is bounded above by $-(1-\eta)\ln (\beta)/(\alpha-1)\ln x$, uniformly over $\ln x\geq x_0$, for suitably large $x_0=x_0(\eta)$. Moreover, it holds that, $c h_n^\delta\beta^{-h_n^\delta}=o(n^{-2})=o(1/\ln x)$, uniformly over $\ln x\leq n^2$, and this completes the proof.
\end{proof}

Finally for this section, we establish that the same distributional tail behaviour for the random variables
\begin{equation}\label{tildeti}
\tilde{t}_i:=\sum_{j=1}^{\tilde{Z}_i-1}\mathbf{1}_{\{h(\mathcal{T}_{ij})\geq h_n\}}\sum_{m=\Delta_i}^{\Delta_{i,(\ln n)^{1+\gamma}}}\mathbf{1}_{\{X_m\in \mathcal{T}_{ij}(x_{ij})\}},
\end{equation}
where $\Delta_{i,(\ln n)^{1+\gamma}}$ is the first time after $\Delta_i$ that the process $X$ hits a backbone vertex outside of the interval $\{\rho_{i-(\ln n)^{1+\gamma}},\dots,\rho_{i+(\ln n)^{1+\gamma}}\}$. Given the backtracking result of Lemma \ref{backtrack}, with high probability it is the case that $\tilde{t}_i$ will be identical to the $t_i$ for all relevant indices $i$. However, the advantage of the sequence $(\tilde{t}_i)$ over $(t_i)$ is that, similarly to the sequence of random variables $(\tilde{T}_x)$ introduced for the directed trap model at (\ref{tildetxdef}), at least when the traps are suitably well-spaced, it is possible to decouple the elements of $(\tilde{t}_i)$ in such a way as to be able to usefully compare them with an independent sequence.

{\lem \label{ttail} Let $\alpha\in(1,2]$ and $\varepsilon>0$. There exist constants $n_0$ and $c_{\beta,\varepsilon}$ such that, for any $i\geq 0$, $n\geq n_0$ and $x$ satisfying $c_{\beta,\varepsilon}h_n \leq \ln x\leq n^2$,
\[ \frac{(1-\varepsilon)\ln \beta}{(\alpha-1)\ln x}\leq \mathbb{P}_\rho\left(\tilde{t}_i\geq x\right)\leq  \frac{(1+\varepsilon)\ln \beta}{(\alpha-1)\ln x}.\]}
\begin{proof} If the process $X$ does not hit $\rho_{i-1-(\ln n)^{1+\gamma}}$ again after having hit $\rho_i$, and does not hit $\rho_{i}$ again after having hit $\rho_{i+1+(\ln n)^{1+\gamma}}$, then $\tilde{t}_i$ is equal to $t_i$. Hence,
\[\mathbb{P}_\rho\left(t_i\neq \tilde{t}_i\right)\leq2\mathbf{E}\left(
P_{\rho_{1+(\ln n)^{1+\gamma}}}^{\mathcal{T}^*}(\tau_{\rho}<\infty)\right).\]
An elementary calculation for the biased random walk on a line shows that the right-hand side here is equal to $\beta^{1-(\ln n)^{1+\gamma}}=o(n^{-2})$. Applying this fact, it is easy to deduce the result from Lemmas \ref{deeptime} and \ref{deeptimelower}.
\end{proof}

\subsection{Proof of main result for critical Galton-Watson trees}

The purpose of this section is to complete the proof of our main results for biased random walks on critical Galton-Watson trees (Theorem \ref{dngrowth} and Corollary \ref{treecor}).

\begin{proof}[Proof of Theorem \ref{dngrowth}] We start the proof by claiming that the conclusion of the lemma holds when the hitting time sequence $(\Delta_{n})_{n\geq0}$ is replaced by $(\sum_{i=0}^{n-1}\tilde{t}_i)_{n\geq 0}$. By imitating the proof of Lemma \ref{tildetxsums} with Lemma \ref{separation} in place of Lemma \ref{prelim1}, to verify that this is indeed the case, it will be enough to prove the same result for
$(\sum_{i=0}^{n-1}\tilde{t}_i')_{n\geq 0}$, where $(\tilde{t}_i')_{i\geq 0}$ in an independent sequence such that $\tilde{t}_i'\sim \tilde{t}_{1+(\ln n)^{1+\gamma}}$ for each $i$. (Note that, because the elements of the sequence $(\tilde{t}_i)_{i\geq 0}$ are only identically-distributed for $i\geq 1+(\ln n)^{1+\gamma}$, we do not take $\tilde{t}_i'\sim \tilde{t}_i$ for each $i$. By applying the second part of Lemma \ref{separation}, which shows that with high probability there will be no big leaves in the interval close to $\rho$, it is easy to adapt the argument of Lemma \ref{tildetxsums} to overcome this issue.) Since the tail asymptotics of Lemma \ref{ttail} mean that the relevant functional scaling limit for $(\sum_{i=0}^{n-1}\tilde{t}'_i)_{n\geq 0}$ is an immediate application of Theorem \ref{thm:Kasa-ext} (with $h_1(n)=\ln n$ and $h_2(n)=n^{-1}$), our claim holds as desired.

Now, fix $T\in(0,\infty)$. By Lemmas \ref{backtrack}, \ref{smallleaf} and \ref{311}, with probability converging to one we have that, for every $t\in [0,T]$,
\[\sum_{i=0}^{nt-1-(\ln n)^{1+\gamma}}\tilde{t}_i\leq \Delta_n\leq \sum_{i=0}^{nt-1}\tilde{t}_i+2\beta^{2h_n}.\]
By repeating the proof of Theorem \ref{onedhit} exactly with the particular choice $L(x):=\log_+x$, this, in conjunction with the conclusion of the previous paragraph,  yields the result.
\end{proof}

\begin{proof}[Proof of Corollary \ref{treecor}] Since the proof is identical to that of Corollary \ref{oned}, with $\bar{F}(x)$ being taken to be a distribution function that is asymptotically equivalent to $\ln \beta/(\alpha-1)\ln x$, we omit it.
\end{proof}

\subsection{Growth rate of quenched mean hitting times}

The purpose of this section is to compare the growth rate of $E^{\mathcal{T}^*}_\rho \Delta_n$, that is, the quenched expectations of the hitting times $\Delta_n$, with the growth rate of $\Delta_n$ that was established in the previous section. Interestingly, in the result corresponding to Theorem \ref{dngrowth} (see Theorem \ref{edngrowth} below), an extra factor of $\alpha$ appears, meaning that the sequence of quenched expectations grows more quickly than the hitting times themselves. This is primarily due to the fact that the quenched expectation $E^{\mathcal{T}^*}_\rho \Delta_n$ feels all the big leaves at a particular backbone vertex, whereas the hitting time $\Delta_n$ only feels the big leaves that are deeply visited by $X$. Indeed, the extra $\alpha$ is most easily understood by comparing the following lemma, which describes the height of the biggest leaf at a particular backbone vertex, with Lemma \ref{hseen}, which concerns only deeply visited big leaves.

{\lem \label{fullheight} Let $\alpha\in(1,2]$. For any $i\geq 0$,
\[\mathbf{P}\left(\max_{j=1,\dots,\tilde{Z}_i-1}h(\mathcal{T}_{ij})\geq x\right)\sim \alpha q_x^{\alpha-1}L\left(q_{x}\right),\]
as $x\rightarrow\infty$.}
\begin{proof} Conditioning on $\tilde{Z}_i$, we obtain
\begin{eqnarray*}
\mathbf{P}\left(\max_{j=1,\dots,\tilde{Z}_i-1}h(\mathcal{T}_{ij})< x\right)&=&\mathbf{E}\left(\mathbf{P}\left(h(\mathcal{T})< x\right)^{\tilde{Z}-1}\right)\\
&=&\mathbf{E}\left(Z\left(1-q_x\right)^{{Z}-1}\right)\\
&=&f'\left(1-q_x\right),
\end{eqnarray*}
where we have once again applied the size-biasing of (\ref{sizebiasz}) to obtain the second equality. Since we know from the proof of Lemma \ref{bigtrap} that $f'(1-x)\sim1-\alpha x^{\alpha-1}L(x)$ as $x\rightarrow 0^+$, the proof is complete.
\end{proof}

In studying the quenched expectation of hitting times, we no longer need an argument that is so sophisticated as to consider the time spent in the individual leaves $\mathcal{T}_{ij}$ (which were defined after \eqref{sizebiasz}). Instead, we will be concerned only with understanding the expected length of time the biased random walk $X$ spends inside sets of the form $\mathcal{T}_{i}=\{\rho_i\}\cup(\cup_{j=1,\dots,\tilde{Z}_i-1}\mathcal{T}_{ij})$. To this end, we introduce a stopping time
\[\sigma_i:=\inf\{n\geq 0:X_n\not\in \mathcal{T}_{i}\}.\]
The expected time spent by $X$ inside $\mathcal{T}_{i}$ on a single visit is thus given by $\mathbf{E}_{\rho_i}^{\mathcal{T}^*}\sigma_{i}$. Similarly to (\ref{treeexit}), we have that
\begin{equation}\label{etbn}
\mathbf{E}_{\rho_i}^{\mathcal{T}^*}\sigma_{i}=\left\{\begin{array}{ll}
                                                     1+\frac{1}{\beta^{i-1}(1+\beta)}\sum_{{x,y\in\mathcal{T}_{i},x\sim y}}c(x,y), & \mbox{if $i\geq 1$}, \\
                                                       1+\sum_{{x,y\in\mathcal{T}_{0},x\sim y}}c(x,y), & \mbox{if $i=0$},
                                                    \end{array}\right.
\end{equation}
and this allows us to obtain the following distributional asymptotics.

{\lem \label{etbntail} Let $\alpha\in(1,2]$. For any $i\geq 0$,
\[\mathbf{P}\left(\mathbf{E}_{\rho_i}^{\mathcal{T}^*}\sigma_i \geq x\right)\sim \frac{\alpha\ln \beta}{(\alpha-1)\ln x},\]
as $x\rightarrow\infty$.}
\begin{proof} If $i\geq1$, then from (\ref{etbn}) we are easily able to deduce that
\begin{equation}\label{taub}
\frac{\beta^{h(\mathcal{T}_i)}}{1+\beta}\leq \mathbf{E}_{\rho_i}^{\mathcal{T}^*}\sigma_i\leq 1+\frac{2\beta^{1+h(\mathcal{T}_i)}\#\mathcal{T}_i}{1+\beta},
\end{equation}
where $h(\mathcal{T}_i)=\mathbf{1}_{\{\tilde{Z}_i> 1\}}+\max_{j=1,\dots,\tilde{Z}_i-1}h(\mathcal{T}_{ij})$ is the height of $\mathcal{T}_{i}$. Hence, for any $\eta\in(0,1)$,
\begin{eqnarray*}
\lefteqn{\mathbf{P}\left(\mathbf{E}_{\rho_i}^{\mathcal{T}^*}\sigma_i \geq x\right)}\\
&\leq & \mathbf{P}\left(\beta^{\max_{j=1,\dots,\tilde{Z}_i-1}h(\mathcal{T}_{ij})}\geq (x-1)^{1-\eta}\right) +\mathbf{P}\left(2\beta^2\#\mathcal{T}_{i}\geq (\beta+1) (x-1)^\eta\right)\\
&\leq &  \mathbf{P}\left(\max_{j=1,\dots,\tilde{Z}_i-1}h(\mathcal{T}_{ij})\geq \frac{(1-\eta)\ln(x-1)}{\ln \beta}\right)+cx^{-\delta\eta},\\
&\sim & \frac{\alpha\ln \beta}{(\alpha-1)(1-\eta)\ln x},
\end{eqnarray*}
where we have applied (\ref{vertexbound}) to deduce the second inequality for suitable constants $c$ and $\delta>0$, and Lemma \ref{fullheight} and (\ref{probdecay}) to obtain the asymptotic equivalence. Since (\ref{taub}) in conjunction with Lemma \ref{fullheight} and (\ref{probdecay}) also implies that
\[\mathbf{P}\left(\mathbf{E}_{\rho_i}^{\mathcal{T}^*}\sigma_i \geq x\right)\geq \mathbf{P}\left(\max_{j=1,\dots,\tilde{Z}_i-1}h(\mathcal{T}_{ij})\geq \frac{\ln x +\ln (1+\beta)}{\ln \beta}\right)\sim \frac{\alpha\ln \beta}{(\alpha-1)\ln x},\]
the result follows in this case. The argument for $i=0$ is similar.
\end{proof}

We are now ready to prove the main result of this section.

{\thm \label{edngrowth} Let $\alpha\in(1,2]$. As $n\rightarrow\infty$, the laws of the processes
\[\left(\frac{(\alpha-1)\ln_+ E^{\mathcal{T}^*}_\rho \Delta_{nt}}{n \alpha \ln \beta}\right)_{t\geq 0}\]
under $\mathbf{P}$ converge weakly with respect to the Skorohod $J_1$ topology on $D([0,\infty),\mathbb{R})$ to the law of $(m(t))_{t\geq 0}$.}
\begin{proof} The embedded random walk on the backbone $Y$ visits each site $\rho_i$, $i\geq 1$, a geometric parameter $(\beta-1)/(\beta+1)$ number of times in total and $\rho=\rho_0$ a geometric parameter $(\beta-1)/\beta$ number of times. Moreover, before visiting $\rho_n$, $Y$ has to visit each element of $\{\rho_0,\dots,\rho_{n-1}\}$ at least once. This and the definition of $\Delta_n$ implies that
\begin{equation}\label{qq}
\sum_{i=0}^{n-1}\mathbf{E}_{\rho_i}^{\mathcal{T}^*}\sigma_i\leq\mathbf{E}_\rho^{\mathcal{T}^*}\Delta_n\leq \frac{\beta+1}{\beta-1}\sum_{i=0}^{n-1}\mathbf{E}_{\rho_i}^{\mathcal{T}^*}\sigma_i.
\end{equation}
Now, the random variables $\mathbf{E}_{\rho_i}^{\mathcal{T}^*}\sigma_i$ in these sums are independent and have slowly varying tails, as described by Lemma \ref{etbntail}. Thus the result is a simple consequence of \cite[Theorem 2.1]{Kasa} (or Theorem \ref{thm:Kasa-ext} below).
\end{proof}

{\rem For comparison, recall the directed trap model of Section \ref{onedd}, but, so as to avoid having to consider the time that the biased random walk $X$ spends at negative integers, replace $\mathbb{Z}$ by the half-line $\mathbb{Z}_+$. As in Theorem \ref{onedhit}, we have that $(n^{-1}L(\Delta_{nt}))_{t\geq 0}$ converges in distribution under the annealed law $\mathbb{P}_0$ to $(m(t))_{t\geq 0}$. For the corresponding quenched expectation, similarly to (\ref{qq}), we have that $\sum_{i=0}^{n-1}\tau_i\leq E^\tau_0(\Delta_n)\leq \frac{\beta+1}{\beta-1}\sum_{i=0}^{n-1}\tau_i$. Thus, again applying 
\cite[Theorem 2.1]{Kasa} (or Theorem \ref{thm:Kasa-ext} below),
it is possible to check that $(n^{-1}L(E^\tau_0(\Delta_{nt})))_{t\geq0}$ converges in distribution under $\mathbf{P}$ to $(m(t))_{t\geq 0}$. In particular, in contrast to the critical Galton-Watson tree case, the asymptotic behaviour of $E^\tau_0(\Delta_n)$ and $\Delta_n$ are identical. This is because, although certain big leaves will be avoided by certain realisations of the biased random walker in the tree setting, the geometry of the graph $\mathbb{Z}_+$ forces $X$, when travelling from $0$ and $n$, to visit all the traps in between on every realisation.}

\section{Extremal aging}\label{agingsec}

In this section, we will prove Theorem~\ref{extagingtree} and Theorem~\ref{extaging}, which state that the biased random walk on critical Galton-Watson tree conditioned to survive and the one-dimensional trap model, respectively, experience extremal aging. The phenomenon we describe for these models is similar to what happens in the trapping models considered by Onur Gun in his PhD thesis~\cite{Gun} and to results observed for spin glasses in~\cite{BG}.

\subsection{Extremal aging for the one-dimensional trap model}

We start by considering the one-dimensional trap model introduced in Section \ref{onedsec}, with the goal of this section being to prove Theorem~\ref{extaging}. The reason for proving this result before its counterpart for trees is that the simpler argument it requires will be instructive when it comes to tackling the more challenging tree case in the subsequent section.

Key to proving Theorem \ref{extaging} is establishing that $X$ localises at the closest trap to 0 of a sufficient depth. To describe this precisely, as we do in Lemma \ref{local} below, we first introduce the notation
\[l(u):=\min\left\{x\geq0:\:\tau_x\geq \bar{F}^{-1}\left(\frac{1}{u}\right)\right\}.\]
From the independently and identically-distributed nature of the environment, we readily deduce the following preliminary lemma.

{\lem \label{lalb} For any $0<a<b$, we have
\[\lim_{n\rightarrow\infty}\mathbf{P}\left(l(an)=l(bn)\right)=\frac{a}{b}.\]}

We now establish the relevant localisation result for $X$.

{\lem \label{local} For any $a>0$, we have
\[\lim_{n\rightarrow \infty}\mathbb{P}_0\left(X_{\bar{F}^{-1}(1/an)}=l(an)\right)=1.\]}
\begin{proof} Our first aim is to show that $X$ hits $l(an)$ before time ${\bar{F}^{-1}(1/an)}$ with high probability. Clearly, for any $T>0$, we have that
\[\mathbf{P}\left(l(an)>nT\right)=\mathbf{P}\left(\tau_x<\bar{F}^{-1}\left(\frac{1}{an}\right):\:x=0,1,\dots,nT\right) = \left(1-\frac{1}{an}\right)^{nT+1}\to e^{-T/a},\]
as $n\rightarrow\infty$. Moreover, by Lemma \ref{lalb}, for $\varepsilon\in(0,1)$ it holds that
\[\mathbf{P}\left(\tau_x\geq \bar{F}^{-1}\left(\frac{1}{an(1-\varepsilon)}\right)\mbox{ for some }x=0,\dots,l(an)-1\right)=\mathbf{P}\left(l(an)\neq l(an(1-\varepsilon))\right)\rightarrow\varepsilon,\]
as $n\rightarrow\infty$. Recalling the notation $\tilde{T}_x$ introduced at~(\ref{tildetxdef}) and applying these two results in conjunction with the bound at (\ref{dnbound}) yields
\[\limsup_{n\rightarrow\infty}\mathbb{P}_0\left(\Delta_{l(an)}> \sum_{x=1}^{nT}\tilde{T}_x\mathbf{1}_{\{\tau_x\leq \bar{F}^{-1}(1/an(1-\varepsilon))\}}+\bar{F}^{-1}\left(n^{-1}(\ln n)^{1/2}\right)\right)\leq \varepsilon+e^{-T/a}.\]
We know that $\bar{F}^{-1}(n^{-1}(\ln n)^{1/2})\leq \bar{F}^{-1}(1/\varepsilon n)$ for large enough $n$, Markov's inequality thus implies
\begin{eqnarray*}
\lefteqn{\limsup_{n\rightarrow\infty}\mathbb{P}_0\left(\Delta_{l(an)}>\bar{F}^{-1}\left(\frac{1}{an}\right)\right)}\\
&\leq & \limsup_{n\rightarrow\infty}\frac{1}{\bar{F}^{-1}(1/an)}\left[
\mathbb{E}_0\left(\sum_{x=1}^{nT}\tilde{T}_x\mathbf{1}_{\{\tau_x\leq \bar{F}^{-1}(1/an(1-\varepsilon))\}}\right)+ \bar{F}^{-1}\left(\frac{1}{\varepsilon n}\right)\right]+\varepsilon+ e^{-T/a}\\
&\leq&\limsup_{n\rightarrow\infty}\frac{1}{\bar{F}^{-1}(1/an)}\left[c_{\beta}nT
\mathbf{E}\left(\tau_0\mathbf{1}_{\{\tau_0\leq \bar{F}^{-1}(1/an(1-\varepsilon))\}}\right)+ \bar{F}^{-1}\left(\frac{1}{\varepsilon n}\right)\right]+\varepsilon+ e^{-T/a},
\end{eqnarray*}
where $c_\beta$ is a constant depending only on $\beta$. By proceeding as in the proof of Lemma \ref{prelim3} with $g(n)$ replaced by $\bar{F}^{-1}(1/an(1-\varepsilon))$, it is possible to check that
\[\mathbf{E}\left(\tau_0\mathbf{1}_{\{\tau_0\leq \bar{F}^{-1}(1/an(1-\varepsilon))\}}\right)\leq \frac{c_1 \bar{F}^{-1}(1/an(1-\varepsilon))}{an(1-\varepsilon)},\]
and so
\[\limsup_{n\rightarrow\infty}\mathbb{P}_0\left(\Delta_{l(an)}>\bar{F}^{-1}\left(\frac{1}{an}\right)\right)
\leq
\limsup_{n\rightarrow\infty}\frac{c_2T\bar{F}^{-1}(1/an(1-\varepsilon))+
a\bar{F}^{-1}(1/\varepsilon n)}{a(1-\varepsilon)\bar{F}^{-1}(1/an)}+\varepsilon+ e^{-T/a}.\]
Were $\limsup_{n\rightarrow\infty} \bar{F}^{-1}(1/an(1-\varepsilon))/\bar{F}^{-1}(1/an)>0$, then there would exist a subsequence $(n_i)_{i\geq 0}$ and constant $c>0$ such that $\bar{F}^{-1}(1/an_i(1-\varepsilon))\geq c \bar{F}^{-1}(1/an_i)$. Applying the decreasing function $\bar{F}$ to both sides and then the slowly-varying property (\ref{svartail}) yields that $a\leq a(1-\varepsilon)$, which is clearly a contradiction. Hence $\lim_{n\rightarrow\infty} \bar{F}^{-1}(1/an(1-\varepsilon))/\bar{F}^{-1}(1/an)=0$. Similarly, one has that $\lim_{n\rightarrow\infty} \bar{F}^{-1}(1/\varepsilon n)/\bar{F}^{-1}(1/an)=0$ for any $\varepsilon<a$. Thus, letting $T\rightarrow \infty$ and $\varepsilon\rightarrow 0$, the above estimate yields
\[\lim_{n\rightarrow\infty}\mathbb{P}_0\left(\Delta_{l(an)}>\bar{F}^{-1}\left(\frac{1}{an}\right)\right)=0,\]
as desired.

Now, if $X_{\bar{F}^{-1}(1/an)}\neq l(an)$, then either $X$ does not hit $l(an)$ before time $\bar{F}^{-1}(1/an)$, or it does hit $l(an)$ and spends less time than $\bar{F}^{-1}(1/an)$ there before moving to any other vertex. By the conclusion of the previous paragraph, the former event has probability 0 asymptotically, and so
\[\limsup_{n\rightarrow\infty}\mathbb{P}_0\left(X_{\bar{F}^{-1}(1/an)}\neq l(an)\right)\leq \limsup_{n\rightarrow\infty}\mathbf{E}\left({P}^\tau_{l(an)}\left(\inf\{t:X_t\neq l(an)\}\leq \bar{F}^{-1}(1/an)\right)\right).\]
Since $\inf\{t:X_t\neq l(an)\}$ is exponential with mean $\tau_{l(an)}$ under ${P}^\tau_{l(an)}$, for any $\varepsilon>0$ the right-hand side here is bounded above by
\[\limsup_{n\rightarrow\infty}\mathbf{E}\left(1\wedge \frac{\bar{F}^{-1}(1/an)}{\tau_{l(an)}}\right)\leq \limsup_{n\rightarrow\infty}\left[
\mathbf{P}\left(\tau_{l(an)}<\bar{F}^{-1}(1/an(1+\varepsilon)) \right)+
\frac{\bar{F}^{-1}(1/an)}{{\bar{F}^{-1}(1/an(1+\varepsilon))}}\right].\]
The probability in the previous expression is equal to $\mathbf{P}(l(an)\neq l(an(1+\varepsilon)))$, and, by Lemma \ref{lalb}, this is asymptotically bounded above by $\varepsilon$. Similarly to an observation made in the previous paragraph, we also have that $\lim_{n\rightarrow\infty} \bar{F}^{-1}(1/an)/\bar{F}^{-1}(1/an(1+\varepsilon))=0$, and thus we have established
\[\limsup_{n\rightarrow\infty}\mathbb{P}_0\left(X_{\bar{F}^{-1}(1/an)}\neq l(an)\right)\leq \varepsilon.\]
Since $\varepsilon$ was arbitrary, this completes the proof.
\end{proof}

Combining Lemmas \ref{lalb} and \ref{local}, we readily obtain Theorem \ref{extaging}.

\subsection{Extremal aging for the critical Galton-Watson tree model}

We now return to the setting of Section \ref{treeint}, so as to prove Theorem~\ref{extagingtree}. Similarly to the strategy of the previous section, we will show that the biased random walk on a critical Galton-Watson tree localises in the first suitably big leaf it visits deeply. To describe this, we introduce the notation:
\[l(x):=\min\left\{i\geq0:\:\max_{j\in V_i}h(\mathcal{T}_{ij})\geq x/\ln \beta\right\}.\]
Whilst the form of the following lemma is similar that of Lemma \ref{lalb}, we note that its proof is more involved. This is because, unlike the holding time means $\tau_x$ used to define $l$ there, the random variables $\max_{j\in V_i}h(\mathcal{T}_{ij})$ are not environment measurable or independent.

{\lem \label{seenheights} Let $\alpha\in (1,2]$. For any $0<a<b$, we have
\[\lim_{n\rightarrow\infty}\mathbb{P}_{\rho}\left(l(an)=l(bn)\right)=\frac{a}{b}.\]}
\begin{proof} First, define
\[\tilde{V}_i:=\left\{j\in B_i:\:\tau_{x_{ij}}<\Delta_{i,(\ln n)^{(1+\gamma)}}\right\},\]
to be the set of big leaves visited by $X$ before the stopping time $\Delta_{i,(\ln n)^{(1+\gamma)}}$ that was introduced at (\ref{tildeti}). Set $\tilde{H}_i:=\max_{j\in \tilde{V}_i}h(\mathcal{T}_{ij})$ if $\tilde{V}_i\neq \emptyset$,
and $\tilde{H}_i=0$ otherwise; observe that if $\tilde{H}_i>0$, then it is necessarily also the case that $\tilde{H}_i\geq h_n$. Moreover, for $T\in(0,\infty)$ and $\varepsilon\in(0,1)$, let
\[\mathcal{E}_1(n):=\left\{\sum_{i=m}^{m+n^\varepsilon}\mathbf{1}_{\{N_n(i)\geq 1\}}\leq 1:\:m=0,1,\dots,Tn-n^\varepsilon
\right\}\cup\left\{\sum_{i=0}^{n^\varepsilon}\mathbf{1}_{\{N_n(i)\geq 1\}}=0\right\}.\]
By proceeding as in the proof of Lemma \ref{tildetxsums}, it is possible to show that, under $\mathbb{P}_\rho$, the random variables $(\tilde{H}_i,N_n(i))_{i=0}^{nT}$ conditional on $\mathcal{E}_1(n)$ have the same joint distribution as $(\tilde{H}'_i,N'_n(i))_{i=0}^{nT}$ conditional on $\mathcal{E}_1'(n)$, where $(\tilde{H}'_i,N'_n(i))_{i\geq 0}$ are independent copies of the pair of random variables $(\tilde{H}_{1+(\ln n)^{1+\gamma}},N_n(1+(\ln n)^{1+\gamma}))$ and $\mathcal{E}_1'(n)$ is defined analogously to $\mathcal{E}_1(n)$ with the $N_n(i)$s replaced by $N_n'(i)$s. Consequently, if we set
\[\tilde{l}(x):=\min\left\{i\geq0:\:\tilde{H}_i\geq x/\ln \beta\right\},\]
and define $\tilde{l}'(x)$ similarly from the random variables $\tilde{H}_i'$, then
\begin{eqnarray}
\lefteqn{\left|\mathbb{P}_\rho\left(\tilde{l}(an)= \tilde{l}(bn)\right)-
\mathbb{P}_\rho\left(\tilde{l}'(an)= \tilde{l}'(bn)\right)\right|}\nonumber\\
&\leq& \left|\mathbb{P}_\rho\left(\tilde{l}(an)= \tilde{l}(bn),\:\tilde{l}(an)\leq nT,\:\mathcal{E}_1(n)\right)-
\mathbb{P}_\rho\left(\tilde{l}'(an)= \tilde{l}'(bn),\:\tilde{l}'(an)\leq nT,\:\mathcal{E}_1'(n)\right)\right|\nonumber\\
&&+\left|\mathbb{P}_\rho\left(\tilde{l}(an)> nT,\:\mathcal{E}_1(n)\right)-\mathbb{P}_\rho\left(\tilde{l}'(an)> nT,\:\mathcal{E}_1'(n)\right)\right|\nonumber\\
&&+2\mathbb{P}_\rho\left(\tilde{l}'(an)> nT,\:\mathcal{E}_1'(n)\right)+2\mathbb{P}_\rho(\mathcal{E}_1(n)^c)\nonumber\\
&\leq &2\mathbb{P}_\rho\left(\tilde{l}'(an)> nT\right)+2\mathbf{P}(\mathcal{E}_1(n)^c),\label{upes}
\end{eqnarray}
where we have applied the fact that $\{\tilde{l}(an)= \tilde{l}(bn),\:\tilde{l}(an)\leq nT\}$ and
$\{\tilde{l}(an)> nT\}$ are both $(\tilde{H}_i)_{i=0}^{nT}$ measurable events. Now, similarly to the observation made in the proof of Lemma \ref{ttail}, if the process $X$ does not hit $\rho_{i-1-(\ln n)^{1+\gamma}}$ again after having hit $\rho_i$, and does not hit $\rho_{i}$ again after having hit $\rho_{i+1+(\ln n)^{1+\gamma}}$ -- an event which has probability greater than $1-o(n^{-2})$ uniformly in $i$, then $\tilde{H}_i$ is equal to $\max_{j\in V_i}h(\mathcal{T}_{ij})$. Hence, applying (\ref{probdecay}) and Lemma \ref{hseen}, we obtain that, for any $x,\varepsilon>0$,
\[\left|\mathbb{P}_\rho\left(\tilde{H}_i\geq xn\right)- \frac{1}{(\alpha-1)xn}\right|\leq \frac{\varepsilon}{n}\]
for large $n$ (uniformly in $i$), and clearly the same bound holds when $\tilde{H}_i$ is replaced by $\tilde{H}'_i$. Applying the independence of the random variables $(\tilde{H}'_i)_{i\geq 0}$, it follows that
\[\lim_{n\rightarrow\infty}\mathbb{P}_\rho\left(\tilde{l}'(an)= \tilde{l}'(bn)\right)=\frac{a}{b},\]
and also
\begin{eqnarray*}
\mathbb{P}_\rho\left(\tilde{l}'(an)> nT\right)
&=&\mathbb{P}_\rho\left(\tilde{H}_{1+(\ln n)^{1+\gamma}}<\frac{an}{\ln \beta}
\right)^{nT+1}\\
&\leq& \left(1-\frac{\ln \beta-\varepsilon a}{(\alpha-1)an}\right)^{nT+1}\\&\sim& e^{- T(\ln \beta-\varepsilon a)/(\alpha-1)a}.
\end{eqnarray*}
Combining these results with Lemma \ref{separation}, which implies that $\mathbf{P}(\mathcal{E}_1(n)^c)\rightarrow 0$, and the estimate at (\ref{upes}), then letting $T\rightarrow \infty$, yields
\begin{equation}\label{tildel}
\lim_{n\rightarrow\infty}\mathbb{P}_{\rho}\left(\tilde{l}(an)=\tilde{l}(bn)\right)=\frac{a}{b}.
\end{equation}

Now, suppose that $\mathcal{E}_2(n)$ is the event that the embedded random walk on the backbone $Y$ does not backtrack more that $(\ln n)^{1+\gamma}$ before hitting $\rho_{n(T+1)}$ -- by Lemma \ref{backtrack}, $\mathbb{P}_{\rho}(\mathcal{E}_2(n))\rightarrow 1$. Moreover, on the event $\mathcal{E}_2(n)$, we have that $\tilde{H}_i=\max_{j\in V_i}h(\mathcal{T}_{ij})$ for $i\leq n(T+1)-1-h_n^\delta$. In particular, for large enough $n$, if $\mathcal{E}_2(n)$ holds and also $\tilde{l}(an)\leq nT$, then it must be the case that $l(an)=\tilde{l}(an)$. Hence, for large $n$,
\[\mathbb{P}_\rho\left(l(an)\neq \tilde{l}(an)\right)
\leq \mathbb{P}_\rho\left(\tilde{l}(an)> nT\right)+\mathbb{P}_\rho\left(\mathcal{E}_2(n)^c\right).\]
Similarly to above, we have that the first term here can be bounded above by
\[\left|\mathbb{P}_\rho\left(\tilde{l}(an)> nT,\:\mathcal{E}_1(n)\right)-\mathbb{P}_\rho\left(\tilde{l}'(an)> nT,\:\mathcal{E}_1'(n)\right)\right|+\mathbb{P}_\rho\left(\tilde{l}'(an)> nT\right)+\mathbf{P}(\mathcal{E}_1(n)^c),\]
the limsup as $n\rightarrow\infty$ of which can be made arbitrarily small by choosing $T$ suitably large. Hence
\[\lim_{n\rightarrow\infty}\mathbb{P}_\rho\left(l(an)\neq \tilde{l}(an)\right)=0.\]
The lemma follows by applying this in conjunction with \eqref{tildel}.
\end{proof}

Before proceeding to prove the analogue of Lemma \ref{local} in the tree setting -- see Lemma \ref{localtree} below, we prove a preliminary estimate which rules out the possibility that any leaves have heights that are close to any particular level on the appropriate scale.

{\lem \label{dddd} Let $\alpha\in(1,2]$. For any $a, T\in(0,\infty)$,
\[\lim_{\varepsilon\rightarrow 0}\limsup_{n\rightarrow\infty}\mathbf{P}\left(\min_{i=0,1,\dots,nT}\min_{j=1,\dots,\tilde{Z}_i-1}\left|h(\mathcal{T}_{ij}) -an\right|\leq a n\varepsilon\right)=0.\]}
\begin{proof} First observe that
\begin{eqnarray*}
\mathbf{P}\left(\min_{j=1,\dots,\tilde{Z}_i-1}\left|h(\mathcal{T}_{ij}) -an\right|\leq a n\varepsilon\right)
&=&1-\mathbf{E}\left(\left(1-q_{an(1-\varepsilon)}+q_{an(1+\varepsilon)}\right)^{\tilde{Z}-1}\right)\\
&=&1-f'\left(1-q_{an(1-\varepsilon)}+q_{an(1+\varepsilon)}\right)\\
&\sim& \alpha \left(q_{an(1-\varepsilon)}-q_{an(1+\varepsilon)}\right)^{\alpha-1}L\left(q_{an(1-\varepsilon)}-q_{an(1+\varepsilon)}\right),
\end{eqnarray*}
where we again apply \cite[(2.1)]{Slack} to deduce the asymptotic equality. Now, from (\ref{probdecay}), one can check that
\[q_{an(1-\varepsilon)}-q_{an(1+\varepsilon)}\sim\frac{2\varepsilon q_{an}}{\alpha-1},\]
which yields
\[\mathbf{P}\left(\min_{j=1,\dots,\tilde{Z}_i-1}\left|h(\mathcal{T}_{ij}) -an\right|\leq a n\varepsilon\right)\sim \alpha\left(\frac{2\varepsilon q_{an}}{\alpha-1}\right)^{\alpha-1}L\left(\frac{2\varepsilon q_{an}}{\alpha-1}\right)\sim \frac{c_\alpha\varepsilon^{\alpha-1}}{an},\]
where $c_\alpha$ is a constant depending only on $\alpha$. The lemma readily follows.
\end{proof}

{\lem \label{localtree} Let $\alpha\in(1,2]$. For any $a>0$, we have
\[\lim_{n\rightarrow \infty}\mathbb{P}_\rho\left(\pi(X_{e^{an}})=\rho_{l(an)}\right)=1.\]}
\begin{proof} Fix $\varepsilon>0$, and let $i_0,j_0$ be indices such that $x_{i_0j_0}$ is the first entrance to a big leaf with height greater than or equal to $an(1+\varepsilon)/\ln \beta$ visited by $X$ (the relevant terminology was introduced just above~(\ref{rrr})). Note that this implies $i_0=l(an(1+\varepsilon))$. Moreover, if $\{l(an(1+\varepsilon))\leq nT\}$ holds and $n$ is suitably large, then $\tau_{x_{i_0j_0}}\leq \Delta_{n(T+1)}$. In particular, for large $n$,
\begin{eqnarray*}
\mathbb{P}_\rho\left(\tau_{x_{i_0j_0}}> \Delta_{n(T+1)}\right)
&\leq& \mathbb{P}_\rho\left(l(an(1+\varepsilon))>nT\right)\\
&\leq& \mathbb{P}_\rho\left(l(an(1+\varepsilon))\neq \tilde{l}(an(1+\varepsilon))\right)+\mathbb{P}_\rho\left(\tilde{l}(an(1+\varepsilon))>nT\right),
\end{eqnarray*}
where $\tilde{l}(an)$ was defined in the proof of Lemma \ref{seenheights}, and the upper bound here converges to 0 as $n$ and then $T$ tend to infinity. Consequently, by applying Lemmas \ref{backtrack}, \ref{smallleaf}, \ref{311} similarly to the proof of Theorem \ref{dngrowth}, as well as Lemma \ref{dddd}, we obtain that if
\[\Theta:=\sum_{i=0}^{n(T+1)-1}\sum_{j= 1}^{\tilde{Z}_i-1}\mathbf{1}_{\{h_n\leq h(\mathcal{T}_{ij})\leq an(1-\varepsilon)/\ln \beta\}}\sum_{m=\Delta_i}^{\Delta_{i,(\ln n)^{1+\gamma}}}\mathbf{1}_{\{X_m\in \mathcal{T}_{ij}(x_{ij})\}},\]
then
\[\lim_{\varepsilon\rightarrow 0}\limsup_{T\rightarrow\infty}\limsup_{n\rightarrow\infty}\mathbb{P}_\rho\left(\tau_{x_{i_0j_0}}>\Theta+2\beta^{2h_n}\right)=0.\]
Now, by proceeding similarly to the proof of Lemma \ref{smallleaf}, we have that
\begin{eqnarray*}
E_{\rho}^{\mathcal{T}^*}\Theta&\leq &\frac{\beta}{\beta-1} \sum_{i=0}^{n(T+1)-1} \tilde{Z}_i\left(1+2e^{an(1-\varepsilon)}{\#\mathcal{T}_{ij}}\right)\hspace{150pt}\\
&\leq &\frac{2\beta n(T+1)e^{an(1-\varepsilon)}}{\beta-1}\max_{i=0,\dots,n(T+1)-1}\tilde{Z}_i\left[1+\max_{j=1,\dots,\tilde{Z}_i-1}\#\mathcal{T}_{ij}\right].
\end{eqnarray*}
Combining these observations yields
\begin{eqnarray*}
\lefteqn{\limsup_{\varepsilon\rightarrow 0}\limsup_{n\rightarrow\infty}\mathbb{P}_\rho\left(\tau_{x_{i_0j_0}}>e^{an}\right)}\\
&\leq & \limsup_{\varepsilon\rightarrow 0}\limsup_{T\rightarrow\infty} \limsup_{n\rightarrow\infty}\mathbb{P}_\rho\left(\Theta>2^{-1}e^{an}\right)\\
&\leq& \limsup_{\varepsilon\rightarrow 0}\limsup_{T\rightarrow\infty} \limsup_{n\rightarrow\infty}\left[\mathbb{P}_\rho\left(\Theta>n E_{\rho}^{\mathcal{T}^*}\Theta\right)+\mathbb{P}_\rho\left(E_{\rho}^{\mathcal{T}^*}\Theta>2^{-1}n^{-1}e^{an}\right)\right]\\
&\leq &  \limsup_{\varepsilon\rightarrow 0}\limsup_{T\rightarrow\infty} \limsup_{n\rightarrow\infty}\mathbf{P}
\left(\max_{i=0,\dots,n(T+1)-1}\tilde{Z}_i\left[1+\max_{j=1,\dots,\tilde{Z}_i-1}\#\mathcal{T}_{ij}\right]>\frac{(\beta-1)e^{an\varepsilon}}{4\beta n^2(T+1)}\right)\\
&\leq & \limsup_{\varepsilon\rightarrow 0} \limsup_{T\rightarrow\infty} \limsup_{n\rightarrow\infty}n(T+1)\mathbf{P}
\left(\tilde{Z}_i\left[1+\max_{j=1,\dots,\tilde{Z}_i-1}\#\mathcal{T}_{ij}\right]>\frac{(\beta-1)e^{an\varepsilon}}{2\beta n^2(T+1)}\right).
\end{eqnarray*}
Applying (\ref{vertexbound}), we thus obtain
\[\lim_{\varepsilon\rightarrow 0}\limsup_{n\rightarrow\infty}\mathbb{P}_\rho\left(\tau_{x_{i_0j_0}}>e^{an}\right)=0.\]

The conclusion of the previous paragraph implies that
\begin{eqnarray*}
\lefteqn{\limsup_{n\rightarrow\infty}\mathbb{P}_\rho\left(\pi(X_{e^{an}})\neq \rho_{l(an)}\right)}\\
&\leq & \limsup_{\varepsilon\rightarrow 0}\limsup_{n\rightarrow\infty}\left[\mathbb{P}_\rho\left(\tau_{x_{i_0j_0}}> e^{an}\right)+\mathbb{P}_\rho\left(\tau_{x_{i_0j_0}}\leq e^{an},\:\pi(X_{e^{an}})\neq \rho_{l(an)}\right)\right]\\
&\leq&\limsup_{\varepsilon\rightarrow 0}\limsup_{n\rightarrow\infty}\mathbf{E}\left(P^{\mathcal{T}^*}_{x_{i_0j_0}}\left(\inf\{m:\:\pi(X_m)\neq \rho_{l(an)}\}\leq e^{an}\right)\right)\\
&\leq & \limsup_{\varepsilon\rightarrow 0}\limsup_{n\rightarrow\infty}\mathbf{E}\left(P^{\mathcal{T}^*}_{x_{i_0j_0}}\left(
\tau_{\rho_{l(an)}}\leq e^{an}\right)\right).
\end{eqnarray*}
It is plain to show that
\[P^{\mathcal{T}^*}_{x_{i_0j_0}}\left(\tau_{y_{i_0j_0}}>\tau_{\rho_{l({an})}}\right)\leq c_1\beta^{-h_{n}^{\delta}},\]
for some constant $c_1$ depending only on $\beta$; indeed, this is nothing more than a computation for a biased random walk on $\mathbb{Z}$. Furthermore, another simple calculation for biased random walk on the line yields
\[P_{y_{i_0j_0}}^{\mathcal{T}^*}\left(\tau_{y_{i_0j_0}}^+>\tau_{\rho_{l(an)}}\right)\leq c_2 e^{-a(1+\varepsilon)n},\]
where $\tau_{y_{i_0j_0}}^+$ is the time of the first return to $y_{i_0j_0}$,
so
\[P_{y_{i_0j_0}}^{\mathcal{T}^*}\left(\tau_{\rho_{l(an)}}\leq e^{an}\right)
\leq \sum_{k=0}^{e^{an}/2}P_{y_{i_0j_0}}^{\mathcal{T}^*}\left(\tau_{y_{i_0j_0}}^+\le \tau_{\rho_{l(an)}}\right)^k
P_{y_{i_0j_0}}^{\mathcal{T}^*}\left(\tau_{y_{i_0j_0}}^+> \tau_{\rho_{l(an)}}\right)\leq c_3e^{-a\varepsilon n}.\]
Consequently,
\[\limsup_{n\rightarrow\infty}\mathbb{P}_\rho\left(\pi(X_{e^{an}})\neq \rho_{l(an)}\right)\leq \limsup_{\varepsilon\rightarrow 0} \limsup_{n\rightarrow\infty}\left( c_1\beta^{-h_{n}^{\delta}}+ c_3e^{-a\varepsilon n}\right)=0,\]
which completes the proof.
\end{proof}

Putting Lemma \ref{seenheights} and \ref{localtree} together, we obtain Theorem \ref{extagingtree}.

\section{A limit theorem for sums of independent random variables with slowly varying tail probability}\label{iidsec}

In this section, we derive the limit theorem for sums of independent
random variables with slowly varying tail probability that was applied in the proofs of Lemma \ref{tildetxsums} and Theorem \ref{dngrowth}. The result we prove here is a generalisation of \cite[Theorem 2.1]{Kasa}.

Let $(X_{i,j})_{i,j\in \mathbb N}$ be non-negative random variables such that for each $n\geq 1$, the elements of the collection $(X_{n,j})_{j\in \mathbb N}$ are independent and have common distribution function $F_n$. Moreover, suppose $F$ is a distribution function such that $\bar F(x):=1-F(x)$ is slowly varying and $\bar F(x)>0$ for all $x>0$. Similarly writing $\bar F_n(x):=1-F_n(x)$, the main assumption of this section is that for each $\varepsilon>0$, there exist constants $c_1,c_2$ such that
\begin{equation} \label{epcond}
(1-\varepsilon)\bar F_n(x)\le \bar F(x)\le (1+\varepsilon)\bar F_n(x),\hspace{20pt}\forall x\in[c_1(g_1(n)\vee 1),c_2g_2(n)],
\end{equation}
where $g_i(n):=\bar{F}^{-1}(n^{-1}h_i(n))$, $i=1,2$, with $h_1: \mathbb N\to (0,\infty)$ a non-decreasing, divergent function satisfying $\lim_{n\to\infty}h_1(n)/n=0$, and $h_2: \mathbb N\to [0,\infty)$ a non-increasing function satisfying $\lim_{n\to\infty}h_2(n)=0$. (Note that necessarily $\lim_{n\to\infty}g_i(n)=\infty$ for $i=1,2$.) Defining a function $L$ by setting $L(x):=1/\bar{F}(x)$, we then have the following scaling result for sums of the form
\[S_m^{n}:=\sum_{j=1}^mX_{n,j}.\]

{\thm\label{thm:Kasa-ext} Assume that (\ref{epcond}) holds. As $n\rightarrow\infty$,
\begin{equation}\label{llimit}
\left(\frac{1}{n} L\left(S_{nt}^n\right)\right)_{t\geq0}\rightarrow\left(m(t)\right)_{t\geq 0}
\end{equation}
in distribution with respect to the Skorohod $J_1$ topology on $D([0,\infty),\mathbb{R})$.}

{\rem (i) Note that, similarly to Remark \ref{extrem}, if $\bar{F}_n$ and $\bar{F}$ are not continuous and eventually strictly decreasing, a minor modification to the proof of the above result (cf. Remark \ref{contrem}(ii)) is needed.\\
(ii) The same conclusion holds if on the left-hand side of (\ref{llimit}) we replace $L$ by $L_n(x)=1/\bar F_n(x)$.}
\bigskip

To the end of proving the above result, it is helpful to introduce $(\eta(t))_{t\geq 0}$ to represent a one-sided stable process with index $1/2$, i.e., with L\'evy measure given by $\mu((x,\infty))=x^{-1/2}$ for $x>0$. We will write $F_*(x)=P(\eta(1)\le x)$ for the distribution function of $\eta(1)$. Briefly, the connection with $(m(t))_{t\geq 0}$ is that $m(t)=(\max_{0\le s\le t}\Delta\eta(s))^{1/2}$ (as processes) -- here we recall $\Delta\eta(s)=\eta(s)-\eta(s^-)$.
Moreover, if we set
\begin{eqnarray*}
\eta_{n,i}&=&\eta(i/n)-\eta((i-1)/n),\hspace{20pt}\forall i\geq 1,\\
m_n(t)&=& \left(\max_{i\le nt}\eta_{n,i}\right)^{1/2}\mathbf{1}_{[1/n,\infty)}(t),
\end{eqnarray*}
then $m_n\rightarrow m$ almost-surely in the Skorohod $J_1$ topology. Indeed, since $\eta_n(t):=\eta(\lfloor nt\rfloor/n)\to \eta(t)$ in the Skorohod $J_1$ topology, by the continuous mapping theorem
\begin{equation}\label{eq:conttheq}
\max_{i\le nt}\eta_{n,i}=\max_{0\le s\le t}\Delta\eta_n(s)\to\max_{0\le s\le t}\Delta\eta(s)
\end{equation}
almost-surely as a process in the same topology.

We are now ready to present the key lemma needed to establish Theorem \ref{thm:Kasa-ext}. (This corresponds to \cite[Lemmas 2.2 and 2.3]{Kasa}.) In its statement, we use the notation $\phi_n(x):=F_n^{-1}(F_*(x))$, and we also define $\phi(x):=F^{-1}(F_*(x))$ for its proof.

{\lem\label{thm:keylem} Under \eqref{epcond}, we have the following.\\
(i) For every $\lambda>0$ and $T>0$, as $n\rightarrow\infty$,
\begin{equation}\label{eq:keyle1}
\sup_{0\le x\le T}\left|\frac 1n L\left(\lambda\phi_n\left(n^2x^2\right)\right)-x\right|\rightarrow0.
\end{equation}
(ii) For each $\delta>0$, $T>\delta$, there exist random constants $K_1,K_2>0$ and $n_0$ such that, for every $t\in [\delta,T]$ and $n\ge n_0$,
\begin{equation}\label{eq:keyle2}
K_1\phi_n\left(n^2m_n(t)^2\right)\le \sum_{i\le nt}\phi_n\left(n^2\eta_{n,i}\right)\le K_2\phi_n\left(n^2m_n(t)^2\right),
\end{equation}
almost-surely.}
\begin{proof}
We first give some preliminary computations. Rewriting \eqref{epcond}, we have
\[\frac {F(x)-\varepsilon}{1-\varepsilon}\le F_n(x)\le \frac{F(x)+\varepsilon}{1+\varepsilon},\hspace{20pt}\forall x\in[c_1(g_1(n)\vee 1),c_2g_2(n)].\]
Setting
\[z=F_*^{-1}\left(F(x)\right),\hspace{20pt}\kappa_1(\varepsilon, z)=F_*^{-1}\left(\frac {F_*(z)-\varepsilon}{1-\varepsilon}\right),\hspace{20pt}
\kappa_2(\varepsilon, z)=F_*^{-1}\left(\frac{F_*(z)+\varepsilon}{1+\varepsilon}\right),\]
it follows that
\begin{equation}\label{eq:kdjeqa}
\phi_n(\kappa_1(\varepsilon, z))\le
\phi(z)=x\le \phi_n(\kappa_2(\varepsilon, z)),
\end{equation}
for $F(c_1g_1(n))<F_*(z)=F(x)<F(c_2g_2(n))$. Since
\[z^{-1/2}\sim1-F_*(z)<1-F(c_1g_1(n))=\bar F(c_1g_1(n))\sim \bar F(g_1(n))=n^{-1}h_1(n),\]
and similarly $z^{-1/2}\sim1-F_*(z)>n^{-1}h_2(n)$, for suitably large $n$, the inequality at \eqref{eq:kdjeqa} holds for all $c_1'n^2/h_1(n)^2<z<c_2'n^2/h_2(n)^2$. Now, since $F_*^{-1}(x)\sim (1-x)^{-2}$ for $x\to 1^-$, we have
\begin{eqnarray*}
\kappa_1(\varepsilon, n^2x^2)&=&F_*^{-1}\left(\frac {F_*(n^2x^2)-\varepsilon}{1-\varepsilon}\right)\\
&\sim& \left(1-\frac {F_*(n^2x^2)-\varepsilon}{1-\varepsilon}\right)^{-2}\\
&=&\frac{(1-\varepsilon)^2}{\left(1-F_*(n^2x^2)\right)^2}\\
&\sim& (1-\varepsilon)^2n^2x^2,
\end{eqnarray*}
so that $\kappa_1(\varepsilon, n^2x^2)\ge (1-\varepsilon)^3n^2x^2$ for large $n$. Similarly, $\kappa_2(\varepsilon, n^2x^2)\le (1+\varepsilon)^3n^2x^2$ for large $n$.
Since $\phi_n, h_1$ are non-decreasing and $h_2$ is non-increasing, by \eqref{eq:kdjeqa} we conclude
\begin{equation}\label{eq:kdjeqa22}
\phi_n\left((1-\varepsilon)^3n^2x^2\right)\le \phi\left(n^2x^2\right)\le \phi_n\left((1+\varepsilon)^3n^2x^2\right),\hspace{20pt}\forall n> h_1^{-1}(c_1''/x)\vee h_2^{-1}(c_2''/x).
\end{equation}

Now let us prove (i). By the definition of $\phi_n$ and the fact that $1-F_*(x)\sim x^{-1/2}$, we have, as in the proof of
\cite[Lemma 2.2]{Kasa}, that $\lim_{x\to\infty}L(\lambda \phi(x^2))/x=1$, from which it follows that $\lim_{n\to\infty}L(\lambda \phi(n^2x^2))/n=x$
for $\lambda>0$. Noting that $L(\lambda \phi(n^2x^2))/n$ is monotone in $x$ and the limiting function is continuous, this convergence is
uniform in $x$ on each finite interval.  By \eqref{eq:kdjeqa22},
\[n^{-1}L\left(\lambda\phi_n\left((1-\varepsilon)^3n^2x^2\right)\right)\le n^{-1}L\left(\lambda\phi(n^2x^2)\right)\le
n^{-1}L\left(\lambda\phi_n\left((1+\varepsilon)^3n^2x^2\right)\right),
\]
for $n>h_1^{-1}(c_1''/x)\vee h_2^{-1}(c_2''/x)$, and thus we obtain \eqref{eq:keyle1}.

We next prove (ii). First, by taking $K_1=1$ the lower bound of \eqref{eq:keyle2} is clear. So we will prove the upper bound.
As in the proof of \cite[Lemma 2.3]{Kasa}, noting that $\phi^{-1}$ is slowly varying and
using the representation theorem of \cite[Theorem 1.2]{sen}, we have $\phi^{-1}(x)=c(x)\exp(\int_1^x \varepsilon (t)/tdt)$ where $c(x)\to c>0$ and $\varepsilon (x)\to 0$ as
$x\to \infty$. Thus, $\phi(x)$ may be expressed as
\[\phi(x)=\exp\left(\int_1^{q(x)x}\frac{\tilde\varepsilon (t)}{t}dt\right)\]
where $\tilde\varepsilon (x)\to \infty$ and $q(x)\to 1/c$ as $x\to\infty$. Using this and \eqref{eq:kdjeqa22}, we have,
for all $a>0$ and $M>2$,
\begin{eqnarray*}
n^2\frac{\phi_n(n^2a)}{\phi_n(2n^2a)}&\le& n^2\frac{\phi\left(\frac{n^2a}{(1-\varepsilon)^{3}}\right)}{\phi_n\left(\frac{2n^2a}{(1+\varepsilon)^3}\right)}\\
&=&n^2\exp\left(-\int_{q\left(\frac{n^2a}{(1-\varepsilon)^{3}}\right)\frac{n^2a}{(1-\varepsilon)^{3}}}
^{q\left(\frac{2n^2a}{(1+\varepsilon)^3}\right)\frac{2n^2a}{(1+\varepsilon)^3}}
\frac{\tilde\varepsilon (t)}{t}dt\right)\\
&\le& n^2\exp(-M\log (an^2))\\
&=&a^{-M}n^{-2(M-1)},
\end{eqnarray*}
when $n> h_1^{-1}(c_1'''/\sqrt a)\vee h_2^{-1}(c_2'''/\sqrt a)\vee (c_Ma^{-1/2})$ for some $c_M$ depending on $M$. Thus
\begin{equation}\label{eq:11-28pm1}
\lim_{n\to\infty}n^2\frac{\phi_n(n^2a)}{\phi_n(2n^2a)}=0,\hspace{20pt}\forall a>0.
\end{equation}
Given this, the rest is a minor modification of the proof of \cite[Lemma 2.3]{Kasa}. Let $a=m(\delta)^2/3>0$. Since
$m_n\to m$ almost-surely (as discussed around \eqref{eq:conttheq}), there exists a random $n_1\ge 1$ such that
\begin{equation}\label{eq:mmnappr}
m_n(t)^2\ge m_n(\delta)^2\ge 2a,\hspace{20pt}\forall t\ge \delta,\:n\ge n_1.
\end{equation}
Let $\mathcal{A}_{1,i}=\{\eta_{n,i}<n^{-2}\}$, $\mathcal{A}_{2,i}=\{n^{-2}\le\eta_{n,i}\le a\}$ and $\mathcal{A}_{3,i}=\{a< \eta_{n,i}\}$, and define
\[S^{n,k}_m=\sum_{i\le m}\phi_n(n^2\eta_{n,i})\mathbf{1}_{\mathcal{A}_{k,i}},\hspace{20pt}k=1,2,3.\]
Since $\phi_n$ is non-decreasing, $S^{n,1}_{nt}\le nT\phi_n(1)$ for $t\le T$. Further, for $n\ge a^{-1/2}$, we have that
$n\phi_n(1)/\phi_n(2n^2a)$ $\le n\phi_n(n^2a)/\phi_n(2n^2a)$, which goes to $0$ as $n\to\infty$ by \eqref{eq:11-28pm1}.
Thus, there exists a random $n_2\ge 1$ such that
\begin{equation}\label{eq:3est-1}
S^{n,1}_{nt}\le \phi_n(2n^2a)\le \phi_n(n^2m_n(t)^2),\hspace{20pt}\forall\delta \le t\le T,\:n\ge n_2,
\end{equation}
where the last inequality is due to \eqref{eq:mmnappr}. Next, by \eqref{eq:11-28pm1}, there exists a random $n_3\ge 1$ such that
\[0\le \phi_n(n^2x)\le \phi_n(n^2a)n^2x\le \phi_n(2n^2a)x,\hspace{20pt}\forall n^{-2}\le x\le a, \:n\ge n_3.\]
Thus, for $\delta\le t\le T$ and $n\ge n_3(\omega)$, we have
\begin{equation}\label{eq:3est-2}
S^{n,2}_{nt}\le \phi_n(2n^2a)\sum_{i\le nt}\eta_{n,i}=\phi_n(2n^2a)\eta(\lfloor nt\rfloor/n)
\le \phi_n(n^2m_n(t)^2)\eta(T),\end{equation}
where the last inequality is due to \eqref{eq:mmnappr}. Now, noting that there are only finitely many $t\in [0,T]$ such that
$\Delta\eta(t)>a$, there exists a random $K_3>0$ such that $\sum_{i\le nT}\mathbf{1}_{A_{3,i}}\le K_3$ for large $n$, almost-surely.
Using this and the definition of $m_n(t)$, there exists a random $n_4\ge 1$ such that the following holds:
\begin{equation}\label{eq:3est-3}
S^{n,3}_{nt}\le K_3\max_{i\le nt}\phi_n(n^2\eta_{n,i})=K_3\phi_n(n^2m_n(t)^2),\hspace{20pt}\forall\delta \le t\le T,\:n\ge n_4.
\end{equation}
Combining \eqref{eq:3est-1}, \eqref{eq:3est-2} and \eqref{eq:3est-3}, we obtain
\[\sum_{i\le nt}\phi_n(n^2\eta_{n,i})\le K_2\phi_n(n^2m_n(t)^2),\hspace{20pt}\forall\delta \le t\le T,\:n\ge n_2\vee n_3\vee n_4=:n_0,
\]
where $K_2=1+\eta(T)+K_3$. Thus we have obtained the upper bound of (ii).
\end{proof}

\begin{proof}[Proof of Theorem \ref{thm:Kasa-ext}]
Given the above lemma, the proof of Theorem \ref{thm:Kasa-ext} is basically the same as that of \cite[Theorem 2.1]{Kasa}, and so we only sketch it briefly. Let
\[\zeta_n^{(n)}(t)=\frac 1n L\left(\sum_{i\le nt}\phi_n\left(n^2\eta_{n,i}\right)\right).\]
Then, by definition, $(\zeta_n^{(n)}(t))_{t\ge 0}$ is equal in law to $(\frac 1n L(S_{nt}^n))_{t\ge 0}$.
Further, as discussed around \eqref{eq:conttheq}, $m_n\to m$ almost-surely with respect to the Skorohod $J_1$ topology. So, in order to complete the proof, it suffices to prove the following: for $T>0$,
\begin{equation}\label{eq:keyaswek}
\sup_{0\le t\le T}\left|\zeta_n^{(n)}(t)-m_n(t)\right|\to 0,
\end{equation}
almost-surely. Firstly,
\begin{eqnarray*}
\sup_{\delta\le t\le T}\left|\zeta_n^{(n)}(t)-m_n(t)\right|&\le& \sup_{0\le t\le T}\max_{i=1,2}\left|n^{-1}L\left(K_i\phi_n(n^2m_n(t)^2)\right)-m_n(t)\right|\\
&\le&\sup_{0\le x\le m_n(T)}\max_{i=1,2}\left|n^{-1}L\left(K_i\phi_n(n^2x^2)\right)-x\right|,\\
\end{eqnarray*}
where Lemma \ref{thm:keylem}(ii) is used in the first inequality. This bound converges to 0 as $n\rightarrow\infty$ by Lemma \ref{thm:keylem}(i). Secondly, using the monotonicity of $\zeta_n^{(n)}$ and $m_n$ (and Lemma \ref{thm:keylem} again), we have
\[\lim_{\delta\to 0}\limsup_{n\to\infty}\sup_{0\le t\le\delta}\left|\zeta_n^{(n)}(t)-m_n(t)\right|\le
\lim_{\delta\to 0}\limsup_{n\to\infty}(\zeta_n^{(n)}(\delta)+m_n(\delta))=\lim_{\delta\to 0}2m(\delta)=0,\]
almost-surely. We thus obtain \eqref{eq:keyaswek}.
\end{proof}

\appendix
\section{Appendix: Skorohod topologies}

For the convenience of the reader, we recall here, following \cite{Whittbook}, the definitions of the two topologies on the Skorohod space $D([0,\infty),\mathbb{R})$ that are applied in this article. We start by defining the Skorohod $J_1$ topology on $D([0,T],\mathbb{R})$, where $T\in (0,\infty)$, to be that induced by the metric
\begin{equation}\label{dj1}
d_{J_1}(f,g):=\inf_{\lambda\in \Lambda}\left\{\|f\circ\lambda-g\|\vee\|\lambda-I\|\right\},
\end{equation}
where $\|\cdot\|$ is the uniform norm, $I$ is the identity map on $[0,T]$ and $\Lambda$ is the set of strictly increasing functions mapping $[0,T]$ onto itself. To define the $M_1$ topology on $D([0,T],\mathbb{R})$, we first introduce the notion of the completed graph of a function $f$ in this space by defining
\[\Gamma_f:=\left\{(t,x)\in[0,T]\times \mathbb{R}:\:x=\alpha f(t^-)+(1-\alpha)f(t)\mbox{ for some }\alpha\in[0,1]\right\}.\]
We then say that $u=(u_1(t),u_2(t))_{t\in[0,1]}$ is a parametric representation of $\Gamma_f$ if it is a continuous bijection from $[0,1]$ to $\Gamma_f$ whose first coordinate is non-decreasing, and define a metric on $D([0,T],\mathbb{R})$ by setting
\[d_{M_1}(f,g):=\inf_{u\in\Pi_f,v\in\Pi_g}\left\{\|u_1-v_1\|\vee\|u_2-v_2\|\right\},\]
where $\Pi_f$ (resp. $\Pi_g$) is the set of parametric representations of $f$ (resp. $g$). It is the topology that $d_{M_1}$ induces that is the Skorohod $M_1$ topology on $D([0,T],\mathbb{R})$. Note that $M_1$ is a weaker topology than $J_1$, in the sense that convergence in the latter implies convergence in former, but not vice versa.

To extend the above notions to $D([0,\infty),\mathbb{R})$, we characterise convergence in the Skorohod $J_1$ (or $M_1$) topology on this space by saying $f_n\rightarrow f$ if and only if $f_n\rightarrow f$ with respect to the Skorohod $J_1$ (or $M_1$) topology on $D([0,T],\mathbb{R})$ for every continuity point of $f$. (These topologies can also be described by metrics, see \cite[Section 3]{Whittbook},
for example.) In particular, to establish weak convergence of a random sequence $(X^n)_{n\geq 1}$ to $X$ with respect to the Skorohod $J_1$ (or $M_1$) topology on $D([0,\infty),\mathbb{R})$, we require that $(X^n)_{n\geq 1}$ converges weakly to $X$ with respect to the Skorohod $J_1$ (or $M_1$) topology on $D([0,T],\mathbb{R})$ for every time $T$ at which $X$ is almost-surely continuous. Note that, since we only ever consider the limits $(m(t))_{t\geq 0}$ and $(m^{-1}(t))_{t\geq 0}$, which are both continuous at each fixed $T$ with probability 1, in our setting we are always required to check that the relevant weak convergence of processes holds in $D([0,T],\mathbb{R})$ for every time $T$.

\section*{Acknowledgements}

Part of this work was completed while D.C. was undertaking a three month JSPS Postdoctoral Fellowship at the Research Institute for Mathematical Sciences, Kyoto University, during which time he was most generously hosted by T.K. A.F.~would like to thank G\'erard Ben Arous for suggesting that these models would exhibit extremal aging.

\providecommand{\bysame}{\leavevmode\hbox to3em{\hrulefill}\thinspace}
\providecommand{\MR}{\relax\ifhmode\unskip\space\fi MR }
\providecommand{\MRhref}[2]{%
  \href{http://www.ams.org/mathscinet-getitem?mr=#1}{#2}
}
\providecommand{\href}[2]{#2}


\begin{thebibliography}{10}

\bibitem{BarKum}
M.~T. Barlow and T.~Kumagai, \emph{Random walk on the incipient infinite
  cluster on trees}, Illinois J. Math. \textbf{50} (2006), no.~1-4, 33--65
  (electronic).

\bibitem{BD}
M.~Barma and D.~Dhar, \emph{Directed diffusion in a percolation network}, J.
  Phys. C: Solid State Phys. \textbf{16} (1983), 1451--1458.

\bibitem{BC}
G.~Ben~Arous and J.~{\v{C}}ern{\'y}, \emph{Dynamics of trap models},
  Mathematical statistical physics, Elsevier B. V., Amsterdam, 2006,
  pp.~331--394.

\bibitem{bacone}
G.~Ben~Arous and J.~{\v{C}}ern{\'y}, \emph{Scaling limit for trap models on {$\mathbb Z^d$}}, Ann. Probab.
  \textbf{35} (2007), no.~6, 2356--2384.

\bibitem{bactwo}
G.~Ben~Arous and J.~{\v{C}}ern{\'y}, \emph{The arcsine law as a universal aging scheme for trap models},
  Comm. Pure Appl. Math. \textbf{61} (2008), no.~3, 289--329.

\bibitem{AFGH}
G.~Ben~Arous, A.~Fribergh, N.~Gantert, and A.~Hammond, \emph{Biased random
  walks on {G}alton{-W}atson trees with leaves}, Ann. Probab. \textbf{40} (2012), no.~1, 280--338.

\bibitem{BG}
G.~Ben~Arous and O.~Gun, \emph{Universality and extremal aging for dynamics of
  spin glasses on sub-exponential time scales}, Comm. Pure Appl. Math.
  \textbf{65} (2012), no.~1, 77--127.

\bibitem{BH}
G.~Ben~Arous and A.~Hammond, \emph{Randomly biased walks on subcritical trees}, Comm. Pure Appl. Math., to appear.

\bibitem{BRWRRW}
D.~A. Croydon, \emph{Slow movement of a random walk on the range of a random
  walk in the presence of an external field}, Preprint, arXiv:1203.0405.

\bibitem{CroyKum}
D.~A. Croydon and T.~Kumagai, \emph{Random walks on {G}alton-{W}atson trees
  with infinite variance offspring distribution conditioned to survive},
  Electron. J. Probab. \textbf{13} (2008), 1419--1441.

\bibitem{Darling}
D.~A. Darling, \emph{The influence of the maximum term in the addition of
  independent random variables}, Trans. Amer. Math. Soc. \textbf{73} (1952),
  95--107.

\bibitem{LegallDuquesne}
T.~Duquesne and J.-F. Le~Gall, \emph{Probabilistic and fractal aspects of
  {L}\'evy trees}, Probab. Theory Related Fields \textbf{131} (2005), no.~4,
  553--603.

\bibitem{Dwass2}
M.~Dwass, \emph{The total progeny in a branching process and a related random
  walk}, J. Appl. Probability \textbf{6} (1969), 682--686.

\bibitem{esz}
N.~Enriquez, C.~Sabot, and O.~Zindy, \emph{Aging and quenched localization for
  one-dimensional random walks in random environment in the sub-ballistic
  regime}, Bull. Soc. Math. France \textbf{137} (2009), no.~3, 423--452.

\bibitem{ESZ2}
N.~Enriquez, C.~Sabot, and O.~Zindy, \emph{Limit laws for transient random walks in random environment on
  {$\mathbb Z$}}, Ann. Inst. Fourier (Grenoble) \textbf{59} (2009), no.~6,
  2469--2508.

\bibitem{Feller2}
W.~Feller, \emph{An introduction to probability theory and its applications.
  {V}ol. {II}.}, Second edition, John Wiley \& Sons Inc., New York, 1971.

\bibitem{FH}
A.~Fribergh and A.~Hammond, \emph{Phase transition for the speed of the biased
  random walk on the supercritical percolation cluster}, Preprint,
  arXiv:1103.1371.

\bibitem{GK}
J.~Geiger and G.~Kersting, \emph{The {G}alton-{W}atson tree conditioned on its
  height}, Probability theory and mathematical statistics ({V}ilnius, 1998),
  TEV, Vilnius, 1999, pp.~277--286.

\bibitem{GneKol}
B.~V. Gnedenko and A.~N. Kolmogorov, \emph{Limit distributions for sums of
  independent random variables}, Addison-Wesley Publishing Company, Inc.,
  Cambridge, Mass., 1954, Translated and annotated by K. L. Chung. With an
  Appendix by J. L. Doob.

\bibitem{Gun}
O.~Gun, \emph{Universality of transient dynamics and aging for spin glasses},
  ProQuest LLC, Ann Arbor, MI, 2009, Thesis (Ph.D.)--New York University.

\bibitem{hmg}
B.~Haas and G.~Miermont, \emph{The genealogy of self-similar fragmentations
  with negative index as a continuum random tree}, Electron. J. Probab.
  \textbf{9} (2004), no. 4, 57--97 (electronic).

\bibitem{Kasa}
Y.~Kasahara, \emph{A limit theorem for sums of i.i.d. random variables with
  slowly varying tail probability}, J. Math. Kyoto Univ. \textbf{26} (1986),
  no.~3, 437--443.

\bibitem{Kesten}
H.~Kesten, \emph{Sub-diffusive behavior of random walk on a random cluster},
  Ann. Inst. H. Poincar\'e Probab. Statist. \textbf{22} (1986), no.~4,
  425--487.

\bibitem{KozNac}
G.~Kozma and A.~Nachmias, \emph{The {A}lexander-{O}rbach conjecture holds in
  high dimensions}, Invent. Math. \textbf{178} (2009), no.~3, 635--654.

\bibitem{rrt}
J.-F. Le~Gall, \emph{Random real trees}, Ann. Fac. Sci. Toulouse Math. (6)
  \textbf{15} (2006), no.~1, 35--62.

\bibitem{LPW}
D.~A. Levin, Y.~Peres, and E.~L. Wilmer, \emph{Markov chains and mixing times},
  American Mathematical Society, Providence, RI, 2009, With a chapter by James
  G. Propp and David B. Wilson.

\bibitem{LPP}
R.~Lyons, R.~Pemantle, and Y.~Peres, \emph{Biased random walks on
  {G}alton-{W}atson trees}, Probab. Theory Related Fields \textbf{106} (1996),
  no.~2, 249--264.

\bibitem{sen}
E.~Seneta, \emph{Regularly varying functions}, Springer-Verlag, Berlin, 1976,
  Lecture Notes in Mathematics, Vol. 508.

\bibitem{Slack}
R.~S. Slack, \emph{A branching process with mean one and possibly infinite
  variance}, Z. Wahrscheinlichkeitstheorie und Verw. Gebiete \textbf{9} (1968),
  139--145.

\bibitem{Whitt}
W.~Whitt, \emph{Weak convergence of first passage time processes}, J. Appl.
  Probability \textbf{8} (1971), 417--422.

\bibitem{Whittbook}
W.~Whitt, \emph{Stochastic-process limits}, Springer Series in Operations
  Research, Springer-Verlag, New York, 2002, An introduction to
  stochastic-process limits and their application to queues.

\bibitem{Zindy}
O.~Zindy, \emph{Scaling limit and aging for directed trap models}, Markov
  Process. Related Fields \textbf{15} (2009), no.~1, 31--50.
\end{thebibliography}
\end{document}